\documentclass[twoside,a4paper,11pt]{amsart}
\usepackage[twoside,asymmetric,bottom=3.5cm,bindingoffset=0pt,nomarginpar]{geometry}

\usepackage{mathtools}
\usepackage{amsmath}
\allowdisplaybreaks[4]
\usepackage{amsthm}
\usepackage{amssymb}
\usepackage{mathrsfs}
\usepackage[all]{xy}
\usepackage{tikz,tikz-cd}
\usepackage[bottom]{footmisc}
\usepackage{bm}
\usepackage{fancyhdr}%
\usepackage{geometry}
\usepackage{hyperref}
\usepackage{cite}
\usepackage{enumitem}
\usepackage{verbatim}
\usepackage{makecell}
\usepackage{fancyhdr}
\setlist[enumerate]{label=(\arabic*)}

\theoremstyle{plain}

\newtheorem{theorem}{Theorem}[section]

\newtheorem{lemma}[theorem]{Lemma}

\newtheorem{corollary}[theorem]{Corollary}

\newtheorem{proposition}[theorem]{Proposition}

\theoremstyle{definition}

\newtheorem{definition}[theorem]{Definition}

\theoremstyle{remark}
\newtheorem{remark}[theorem]{Remark}

\newcommand{\hU}{\hat{U}}

\newcommand{\R}{\mathbb{R}}
\newcommand{\RR}{\mathbb{R}}

\DeclareSymbolFont{AMSb}{U}{msb}{m}{n}
\DeclareMathSymbol{\Z}{\mathalpha}{AMSb}{"5A}
\DeclareMathSymbol{\D}{\mathalpha}{AMSb}{"44}
\DeclareMathSymbol{\s}{\mathalpha}{AMSb}{"53}

\numberwithin{equation}{section}
\allowdisplaybreaks

\newcommand{\nn}{\mathcal{N}}

\newcommand{\orb}{\mathcal{O}}
\newcommand{\Or}{\mathcal{O}}

\newcommand{\Ker}{{\text{Ker}}}
\newcommand{\rank}{{\text{rank}}}

\newcommand{\diam}{{\text {diam}}}

\newcommand{\supp}{{\text {supp}}}

\newcommand{\iso}{\text{Iso}}

\newcommand{\thi}{\text{thick}}

\newcommand{\eps}{\epsilon}
\newcommand{\rcd}{\text{RCD}}

\newcommand{\Ric}{\text{Ric}}
\newcommand{\mm}{\mathfrak{m}}

\newcommand{\op}[1]{\operatorname{#1}}

\begin{document}
	\title[Rigidity and regularity for almost homogeneous spaces with Ricci curvature bounds] {Rigidity and regularity for almost homogeneous spaces with Ricci curvature bounds}
	
	\author{Xin Qian}
	\address{School of Mathematical Sciences, Fudan University, Shanghai China}
	\curraddr{}
	\email{xqian22@m.fudan.edu.cn}


	\begin{abstract}
		We say that a metric space $X$ is $(\epsilon,G)$-homogeneous if $G\leq\op{Iso}(X)$ is a discrete group of isometries with $\op{diam}(X/G)\leq\epsilon$.\ A sequence of $(\epsilon_i,G_i)$-homogeneous spaces $X_i$ with $\epsilon_i\to0$ is called a sequence of almost homogeneous spaces. 
		
		In this paper we show that the Gromov-Hausdorff limit of a sequence of almost homogeneous RCD$(K,N)$ spaces must be a nilpotent Lie group with $\Ric \geq K$.\ We also obtain a topological rigidity theorem for $(\epsilon,G)$-homogeneous RCD$(K,N)$ spaces, which generalizes a recent result by Wang.\ Indeed, if $X$ is an $(\epsilon,G)$-homogeneous RCD$(K,N)$ space and $G$ is an almost-crystallographic group, then $X/G$ is bi-H\"older to an infranil orbifold.\ Moreover, we study $(\epsilon,G)$-homogeneous spaces in the smooth setting and prove rigidity and $\epsilon$-regularity theorems for Riemannian orbifolds with Einstein metrics and bounded Ricci curvatures respectively.
	\end{abstract}

	\maketitle

	\tableofcontents
	\vspace{-5pt}
	
	\section{Introduction}\label{section-0}
	
	A classical result of Gromov \cite{gromov1978almost} (refined by Ruh \cite{ruh1982almost}) on almost flat manifolds states that for any integer $n\geq2$, there exists $\epsilon(n), C(n)>0$ such that if a closed $n$-manifold $(M,g)$ satisfies $\op{diam}(M,g)\leq\epsilon(n)$ and $|\op{sec}_g|\leq 1$, then $M$ is diffeomorphic to a infranilmanifold $\nn/\Gamma$, where $\nn$ is a simply connected nilpotent Lie group, and $\Gamma$ is a torsion free discrete subgroup of affine group $\nn\rtimes\op{Aut}(\nn)$ such that $[\Gamma:\Gamma\cap \nn]\le C(n)$.\ This topological control fails if one works on manifolds with bounded Ricci curvature, since even the topology of compact Ricci-flat manifold can exhibit considerable complexity.
	
	However, the nilpotent structure still occurs at the level of fundamental group for manifolds with lower bounded Ricci curvature.\ In fact, it was proved by Kapovitch-Wilking \cite{kapovitch2011structure} that there exists $\epsilon(n), C(n)>0$ such that for any closed $n$-manifold $(M,g)$ with $\op{diam}(M,g)\leq\epsilon(n)$ and $\op{Ric}_g\geq-(n-1)$, its fundamental group $\pi_1(M)$ contains a nilpotent subgroup $N$ with index $[\pi_1(M):N]\le C(n)$ and $\op{rank}(N)\leq n$.\ This theorem is now called the generalized Margulis lemma and it has recently been generalized in the non-smooth setting (i.e., for RCD spaces) in \cite{deng2023margulis}.\ It is known that any finitely generated nilpotent group is polycyclic and $\op{rank}(N)$ is defined as the number of $\mathbb Z$ factors in the polycyclic series of $N$.\ In general, following \cite{zamora2024topological}, we can define the rank for any finitely generated group $G$ as the infimum of $\op{rank}(N)$ among all finite index polycyclic subgroups $N\leq G$ (see Definition \ref{def-rank}).\ In particular, for a closed $n$-manifold $(M,g)$ with $\op{diam}(M,g)\leq\epsilon(n)$ and $\op{Ric}_g\geq-(n-1)$, $\op{rank}(\pi_1(M))\leq n$.
	
	By Naber-Zhang's results in \cite{naber2016topology}, if $\pi_1(M)$ attains maximal rank $n$, then the universal cover $\widetilde{M}$ is volume non-collapsing; Rong termed this the manifold $M$ satisfying the bounded covering geometry (see \cite{huang2020collapsed,rong2022collapsed}).\ It was proved in \cite{huang2020collapsed} that the bounded covering geometry will result in $M$ being an infra-nilmanifold, which generalized Gromov-Ruh's theorem on almost flat manifolds.\ Combining the results in \cite{naber2016topology,huang2020collapsed}, we have the following theorem, which was also mentioned in a recent work by Si-Xu \cite{si2024rigidity}.
	
	\begin{theorem}[\citeonline{naber2016topology,huang2020collapsed}]\label{thm-almostflat}
		There is $\epsilon(n)>0, v(n)>0$ such that for any $n$-manifold
		$(M,g)$ with $\op{Ric}_g\geq -(n-1)$ and $\op{diam}(M)\leq\epsilon(n)$, the followings are equivalent:
		\begin{enumerate}
			\item $M$ is diffeomorphic to an infranilmanifold;
			\item $\op{rank}(\pi_1(M))$ is equal to $n$;
			\item $(M,g)$ satisfies $(1,v(n))$-bounded covering geometry, i.e., $\op{vol}(B_1(\tilde x))\geq v(n)$, where $\tilde x$ is a point in the Riemannian universal cover $\widetilde{M}$.
		\end{enumerate}
	\end{theorem}
	
	In \cite{si2024rigidity}, Si-Xu proved that if in addition $(M,g)$ is Einstein in Theorem \ref{thm-almostflat}, then $(M,g)$ must be flat.\ Moreover, based on the results of the two most recent works by Zamora-Zhu \cite{zamora2024topological} and Wang \cite{wang2024non}, which will be  reviewed in Section \ref{preliminaries}, Theorem \ref{thm-almostflat} can actually be extended to the RCD setting as following.
	
	\begin{theorem}[\citeonline{wang2024non,zamora2024topological}]\label{thm-rcd-almostflat}
		For each $K\in\mathbb{R}$ and $N\geq1$, there is $\epsilon=\epsilon(K,N), v=v(K,N)$ such that for any RCD($K,N$) space $(X,d,\mm)$ with $\op{diam}(X)\leq\epsilon$, the followings are equivalent:
		\begin{enumerate}
			\item $X$ is bi-H\"older homeomorphic to an $N$-dimensional infranilmanifold $\mathcal{N}/\Gamma$, where $\mathcal{N}$ is a simply connected nilpotent Lie group endowed with a left invariant metric;
			\item $\op{rank}(\pi_1(X))$ is equal to $N$;
			\item $\mathcal{H}^N(B_1(\tilde x))\geq v$, where $\tilde x$ is a point in the universal cover $\widetilde{X}$.
		\end{enumerate}
	\end{theorem}
	
	In this paper, we aim to study such questions in a more general situation.\ Accordingly, we first introduce the following definition (see \cite{zamora2024limits}).
	
	\begin{definition}\label{def-almosthomogeneous}
		We say a proper geodesic space $X$ is $(\epsilon,G)$-homogeneous if $G\leq\op{Iso}(X)$ is a discrete group of isometries with $\op{diam}(X/G)\leq\epsilon$.\ A sequence of $(\epsilon_i,G_i)$-homogeneous spaces $X_i$ with $\epsilon_i\to0$ is called a sequence of almost homogeneous spaces. 
	\end{definition}
	For example, if $X$ is a compact geodesic space with $\op{diam}(X)\leq\epsilon$, then any normal cover $\hat{X}$ is $(\epsilon,G(\hat{X}))$-homogeneous where $G(\hat{X})$ is the group of deck transformations.\ In particular, the universal cover $\widetilde{X}$ is $(\epsilon,\pi_1(X))$-homogeneous.\ In general, the action of $G$ may not be free.\ Notice that we require discreteness for the group $G$ in Definition \ref{def-almosthomogeneous} and hence, a homogeneous space is not necessarily an $(\epsilon,G)$-homogeneous space.
	
	The goal of this paper is to study the structure of $(\epsilon,G)$-homogeneous spaces with Ricci curvature bounds.\ In the case of lower Ricci curvature bounds, we investigate this in the general non-smooth setting (i.e., the RCD setting).\ However, for two-sided Ricci curvature bounds or Einstein metrics, we are limited to a smooth setting, focusing on Riemannian manifolds and orbifolds.\ Our first main result is to classify the pointed measured GH-limit of a sequence of almost homogeneous RCD($K,N$) spaces.
	\begin{theorem}\label{thm-almost homogeneous rcd}
		Let $K\in\mathbb{R}, N\in[1,\infty)$ and $(X_i,d_i,\mm_i,p_i)$ be a sequence of $(\epsilon_i,G_i)$-homogeneous pointed RCD$(K,N)$ spaces converging to $(X,d,\mm,p)$ in the pmGH-sense with $\epsilon_i\to0$.\ Assume that $X$ is not a point.\ Then the followings hold:
		\begin{enumerate}
			\item $(X,d)$ is isometric to an $n$-dimensional nilpotent Lie group with a left invariant Riemannian metric for some $n\leq N$;
			\item if the actions of $G_i$ are measure-preserving, then $\mm=c\mathcal{H}^n$ for some $c>0$ and $(X,d,\mathcal{H}^n)$ is an RCD$(K,n)$ space.\ In particular, $\op{Ric}_X\geq K$ if $n\geq2$;
			\item if $X_i$ are compact (or equivalently, $G_i$ are finite), then $\mm$ and $\mathcal{H}^n$ are mutually abolutely continuous and $(X,d,\mathcal{H}^n)$ is an RCD$(K,n)$ space.\ In particular, $\op{Ric}_X\geq K$ if $n\geq2$;
			\item if $X$ is compact, then $X$ is isometric to a flat torus $\mathbb{T}^n$.
		\end{enumerate}
	\end{theorem} 
	For a metric measure space $(Y,d,\mm)$, we say an isometric action $g\in\op{Iso}(Y)$ is measure preserving if $g_{\#}\mm=\mm$.\ Notice that when $\mm$ is the Hausdorff measure, any isometric action is measure preserving.\ Also, if we consider a normal covering $(\hat{Y},\hat{d},\hat{\mm})\to(Y,d,\mm)$, then the construction of the lifted measure $\hat{\mm}$ trivially implies that the deck transformations are measure preserving.
	
	Theorem \ref{thm-almost homogeneous rcd} implies that the pGH-limit of a sequence of pointed almost homogeneous manifolds with $\Ric\geq K$ must be isometric to one of the following: a point, $S^1(r)$, $\mathbb{R}$, or a nilpotent Lie group with left invariant Riemannian metric and $\Ric\geq K$.
	
	Based on Gigli's splitting theorem on RCD($0,N$) spaces (see Theorem \ref{gigli}), we also obtain the following result.
	\begin{theorem}\label{thm-almost nonnegative ricci-limit}
		Let $(X_i,d_i,\mm_i,p_i)$ be a sequence of almost homogeneous RCD$(-\delta_i,N)$ spaces converging to $(X,d,\mm,p)$ in the pmGH-sense with $\delta_i\to0$.\ Then $X$ is isometric to $\mathbb{R}^{k}\times\mathbb{T}^{n-k}$, where $\mathbb{T}^{n-k}$ is a flat torus and $0\leq k\leq n\leq N$.
	\end{theorem}
	

	We can also prove the following corollary for non-collapsed RCD spaces.\ Recall that an RCD($K,N$) space $(X,d,\mm)$ is called non-collapsed if $\mm=\mathcal{H}^N$ (see \cite{de2018non}).\ In this case, $N$ must be an integer.
	
	\begin{corollary}\label{cor-homeomorphic to tori}
		Given $v>0,D>0,K\in\mathbb{R}$ and $N\geq1$, there exists $\epsilon=\epsilon(K,N,v,D)$ such that if $(X,d,\mathcal{H}^N)$ is an $(\epsilon,G)$-homogeneous RCD($K,N$) space with $\mathcal{H}^N(X)\geq v$ and $\diam(X)\leq D$, then $X$ is bi-H\"older homeomorphic to a flat torus $\mathbb{T}^N$.\ Moreover, if $X$ is a Riemannian manifold, then $X$ is diffeomorphic to $\mathbb{T}^N$.
	\end{corollary}
	
	More generally, one may consider those $(X,d,\mathcal{H}^N)$ satisfying local $(r,v)$-bounded covering geometry; that is, for any $x\in X$, $\mathcal{H}^N(B_{r/2}(\tilde{x}))\geq v$, where $\tilde{x}$ is a pre-image of $x$ in the (incomplete) universal covering $\widetilde{B_r(x)} \to B_r(x)$.\ Notice that on an RCD$(K,N)$ space, any ball $B_r(x)$ admits a universal cover by Wang's result \cite{wang2024rcd} (see Theorem \ref{jikang}).\ Employing Wang's latest results \cite{wang2024non}, we have the following fibration theorem for almost homogeneous RCD$(K,N)$ spaces satisfying local bounded covering geometry.
	
	\begin{theorem}\label{thm-covering and almost homogeneous}
		Given $v>0,D>0,K\in\mathbb{R}$ and $N\geq1$, there exists $\eps=\eps(K,N,v,D)$ such that if $(X,d,\mathcal{H}^N)$ is an $(\epsilon,G)$-homogeneous RCD($K,N$) space satisfying local $(1,v)$-bounded covering geometry and $\diam(X)\leq D$, then there is a flat torus $\mathbb{T}^k$ with $0\leq k\leq N$ and a continuous fiber bundle map $f: X\to\mathbb{T}^k$, which is also an $\Phi(\eps|K,N,v,D)$-GHA, where $\Phi(\eps|K,N,v,D)\to0$ as $\epsilon\to0$ while $K,N,v$ and $D$ are fixed.\ Moreover, the $f$-fiber with the induced metric is bi-H\"older to an $(N-k)$-dimensional infranilmanifold and in particular, $X$ is homeomorphic to an infranilmanifold.
		
		In addition, if $X$ is a Riemannian manifold, then $f$ is a smooth bundle map and the $f$-fiber is diffeomorphic to an $(N-k)$-dimensional infranilmanifold.\ In particular, $X$ is diffeomorphic to an infranilmanifold
	\end{theorem}
	
	Notice that $\diam(X)\leq\eps$ is equivalent to $X$ being $(\eps,\left\lbrace e\right\rbrace)$-homogeneous.\ In this case, the flat torus $\mathbb{T}^k$ in Theorem \ref{thm-covering and almost homogeneous} will be a point, and thus this theorem can be seen as a generalization of \cite[Theorem A]{huang2020collapsed} and \cite[Theorem A]{wang2024non}, which are the “(3) implies (1)" parts of Theorem \ref{thm-almostflat} and Theorem \ref{thm-rcd-almostflat}, respectively.

	For a proper geodesic space $X$, $G\leq\op{Iso}(X)$ is discrete if and only if its action on $X$ has discrete orbits and is almost free (any isotropy group is finite).\ If $X$ is a manifold, then $X/G$ admits a natural orbifold structure \cite{thurston1997three}.\ Recall that an infranil orbifold is the quotient of a simply connected nilpotent Lie group $\mathcal{N}$ by an almost-crystallographic group (see Definition \ref{def-almost crystal}).\ We can then generalize Theorem \ref{thm-rcd-almostflat} as following.
	
	\begin{theorem}\label{thm-rcd-maximal rank}
		For each $K\in\mathbb{R}$ and $N\geq1$, there is $\epsilon=\epsilon(K,N), v=v(K,N)$ such that for any $(\eps,G)$-homogeneous RCD($K,N$) space $(X,d,\mm)$, the followings are equivalent:
		\begin{enumerate}
			\item $X$ is homeomorphic to $\mathbb{R}^N$;
			\item $X$ is a contractible topological $N$-manifold without boundary;
			\item $\op{rank}(G)$ is equal to $N$;
			\item $X$ is simply connected and $\mathcal{H}^N(B_1(x))\geq v$ for some $x\in X$;
			\item $\pi_1(X)$ is finite and $\mathcal{H}^N(B_1(x))\geq v$ for some $x\in X$.
		\end{enumerate}
		Moreover, if $G$ does not contain a non-trivial finite normal subgroup, then the above conditions are equivalent to the following:
		\begin{enumerate}
			\item[(6)] $X/G$ is bi-H\"older homeomorphic to an $N$-dimensional infranil orbifold $\mathcal{N}/\Gamma$, where $\mathcal{N}$ is a simply connected nilpotent Lie group endowed with a left invariant metric and $G$ is isomorphic to $\Gamma$.
		\end{enumerate}
	\end{theorem}
	If $X\to X/G$ is a covering map, then $G$ is torsion free when $X$ is a contractible manifold.\ Thus, Theorem \ref{thm-rcd-maximal rank} will indeed yield Theorem \ref{thm-rcd-almostflat} (see Remark \ref{cor-rcd}).\ In general, we need to assume that $G$ does not contain a non-trivial finite normal subgroup to identify $G$ as an almost-crystallographic group (see Theorem \ref{thm-almost crystallographic}).
	
	
	The construction of the bi-H\"older homeomorphism in Theorem \ref{thm-rcd-maximal rank} (6) is based on the recent work of Wang \cite{wang2024non} and can be adapted to prove diffeomorphism in the smooth case.\ In fact, the proof of Theorem \ref{thm-rcd-maximal rank} allows us to obtain the following orbifold version of Theorem \ref{thm-almostflat}.\ Basic notions and terminology regarding smooth (Riemannian) orbifolds will be reviewed in Section \ref{subsec:orb}. 
	
	\begin{theorem}\label{thm-orbifold almost flat}
		There is $\epsilon(n)>0, v(n)>0$ such that if $(\mathcal{O},g)$ is a Riemannian $n$-orbifold with $\op{Ric}\geq-(n-1)$ and $\op{diam}(\mathcal{O})\leq\epsilon(n)$, then the followings are equivalent:
		\begin{enumerate}
			\item $|\tilde \orb|$ is homeomorphic to $\mathbb{R}^n$, where $|\tilde \orb|$ is the underlying topological space of the universal orbifold cover $\tilde{\orb}$;
			\item $\op{rank}(\pi_1^{orb}(\mathcal{O}))=n$, where $\pi_1^{orb}(\orb)$ is the orbifold fundamental group of $\orb$;
			\item $\op{vol}(B_1(\tilde x))\geq v(n)$, where $\tilde x$ is a point in the universal orbifold covering $\widetilde{\mathcal{O}}$.
		\end{enumerate}
		Moreover, if $\orb$ is good and $\pi_1^{orb}(\orb)$ does not contain a non-trivial finite normal subgroup, then the above conditions are equivalent to the following:
		\begin{enumerate}
			\item[(4)] $\orb$ is diffeomorphic to an infranil orbifold.
		\end{enumerate}
	\end{theorem}
	Recall that a smooth orbifold is called good if it is the global quotient of a smooth manifold by some dicrete group.\ In the last statement of Theorem \ref{thm-orbifold almost flat}, without assuming $\orb$ is good, we know from Theorem \ref{thm-rcd-maximal rank} that the underlying space $|\orb|$ is homeomorphic to an infranil orbifold; with this assumption, we obtain diffeomorphism.

	For two-sided Ricci curvature bounds or Einstein metrics, there is currently no comprehensive and rigorous synthetic theory on metric measure spaces.\ So we are limited to a smooth setting and the following theorem is about the rigidity case in Theorem \ref{thm-orbifold almost flat} for Einstein orbifolds, which can be seen as an orbifold version of \cite[Theorem 0.2]{si2024rigidity}.
	
	

	\begin{theorem}[Rigidity for Einstein orbifolds]\label{thm-einstein orbifold}
		There is $\epsilon(n)>0, v(n)>0$ such that for any closed Einstein $n$-orbifold $(\mathcal{O},g)$ satisfying $\op{Ric}\equiv\lambda g$ with $\lambda\geq-(n-1)$ and $\op{diam}(\mathcal{O})\leq\epsilon(n)$, the followings are equivalent:
		\begin{enumerate}
			\item $\mathcal{O}$ is a closed flat orbifold;
			\item $\orb$ is diffeomorphic to an infranil orbifold;
			\item $\op{rank}(\pi_1^{orb}(\mathcal{O}))=n$;
			\item $\op{vol}(B_1(\tilde x))\geq v(n)$, where $\tilde x$ is a point in the universal orbifold covering $\widetilde{\mathcal{O}}$.
		\end{enumerate}
	\end{theorem}
	
	Compared with Theorem \ref{thm-orbifold almost flat}, we removed the assumption that $\orb$ is good and $\pi_1^{orb}(\mathcal{O})$ does not contain a non-trivial finite normal subgroup, since we can get a flat metric in the Einstein case.\ Recall that by a result of Thurston \cite{thurston1997three}, a closed flat orbifold must be the quotient orbifold $\mathbb{R}^n/\Gamma$, where $\Gamma$ is a crystallographic group.
	
	It is known that almost flat orbifolds are infranil orbifolds (see \cite[Proposition 1.4]{ding2011restriction}).\ Then it follows from Theorem \ref{thm-einstein orbifold} that almost flat Einstein orbifolds must be flat.
	
	\begin{corollary}
		There exists $\epsilon(n)>0$ such that if $(\orb,g)$ is an Einstein $n$-orbifold satisfying $\op{diam}(\orb)^2 |\sec_g|\leq \epsilon(n)$, then $(\orb,g)$ is flat.
	\end{corollary}
	
	For the case of bounded Ricci curvature, we have the following $\epsilon$-regularity theorem.
	
	\begin{theorem}\label{thm-regularity}
		Given $v>0$ and $p\in(1,\infty)$, there exists $C=C(n,v,p), \epsilon=\epsilon(n,v,p)$ such that if an $(\epsilon,G)$-homogeneous pointed $n$-orbifold $(\mathcal{O},g,x)$ satisfies $|\op{Ric}_g|\leq n-1$ and $\op{vol}(B_1(x))\geq v$, then $\int_{B_1(x)}|Rm|^p\leq C$.\ If, in addition, $(\orb,g)$ is a manifold, then there exists $r_0=r_0(n,v)$ such that $r_h(y)\geq r_0$ for any $y\in B_1(x)$, where $r_h(y)$ is the $C^1$-harmonic radius at $y$.
	\end{theorem}

	For a manifold with bounded Ricci curvature, a lower bound on $C^1$-harmonic radius implies a lower bound on $C^{1,\alpha}$-harmonic radius and a local $L^p$-bound on curvature, due to elliptic estimates.\ For Einstein manifolds, this will further lead to higher-order control on curvature.\ However, under our conditions, the Einstein metric exhibits very strong rigidity and will, in fact, be flat by Theorem \ref{thm-einstein orbifold}.
	
	Moreover, for manifolds with bounded Ricci curvature, we obtain the following result analogous to Theorem \ref{thm-almostflat}.
	
	\begin{theorem}\label{thm-bounded Ricci almostflat}
		Given $p>n/2$, there is $C(n,p), \epsilon(n,p), v(n)$ such that for any $n$-manifold
		$(M,g)$ with $|\op{Ric}_g|\leq n-1$ and $\op{diam}(M)\leq\epsilon(n,p)$, the followings are equivalent:
		\begin{enumerate}
			\item $||Rm||_{L^p(M)}\leq C(n,p)$;
			\item $M$ is diffeomorphic to an $n$-dimensional infra-nilmanifold;
			\item $\op{rank}(\pi_1(M))$ is equal to $n$;
			\item $(M,g)$ satisfies $(1,v(n))$-bounded covering geometry, i.e., $\op{vol}(B_1(\tilde x))\geq v(n)$, where $\tilde x$ is a point in the Riemannian universal cover $\widetilde{M}$.
		\end{enumerate}
	\end{theorem}
	The paper is organized as follows.\ In Section \ref{preliminaries}, we cover the preliminary material.\ In Section \ref{sec-3}, we prove Theorem \ref{thm-almost homogeneous rcd}, Theorem \ref{thm-almost nonnegative ricci-limit} and derive a series of consequences.\ In Section \ref{sec-4}, the topological rigidity theorems (Theorem \ref{thm-rcd-maximal rank} and Theorem \ref{thm-orbifold almost flat}) will be proved using results in \cite{zamora2024topological} and \cite{wang2024non}.\ In Section \ref{sec-5}, we study the rigidity and $\epsilon$-regularity for Einstein orbifolds and orbifolds with bounded Ricci curvature respectively.

	\textbf{Acknowledgements:} The author would like to thank Ruobing Zhang for bringing Jikang Wang's recent paper \cite{wang2024non} to his attention.
	
	\section{Preliminaries}\label{preliminaries}
	In this paper, a metric measure space is a triple $(X,d,\mathfrak{m})$, where $(X,d)$ is a complete, separable and proper metric space and $\mathfrak{m}$ is a locally finite non-negative Borel measure on $X$ with $\supp\mathfrak{m}=X$.\ We will also always assume $(X,d)$ to be geodesic, i.e., any couple of points is joined by a length minimizing geodesic.
	
	Throughout this paper, we assume the reader is familiar with the notion and basic theory of $\rcd (K,N)$ spaces ($N<\infty$).\ We refer the reader to \cite{lott2009ricci,sturm2006geometry,sturm2006geometry2,ambrosio2014metric,ambrosio2014calculus,gigli2015differential} for the relevant notions.\ Notice that a large body of literature has studied the so-called RCD$^*(K,N)$ spaces, which are now known to be equivalent to RCD$(K,N)$ spaces by the work of Cavelletti-Milman \cite{cavalletti2021globalization} and Li \cite{li2024globalization}.\ We refer to \cite{erbar2015equivalence} for an overview of equivalent definitions of $\rcd (K,N)$ spaces.

	\subsection{Gromov--Hausdorff topology}
	
	\begin{definition}
		Let $(X,p),(Y,q)$ be pointed geodesic spaces.\ A map $f:X\to Y$ is called a pointed $\epsilon$-Gromov-Hausdorff approximation (or pointed $\epsilon$-GHA) if 
		\begin{equation}\label{pgh1}
			d(f(p),q)\leq\eps, 
		\end{equation}
		\begin{equation}\label{pgh2}
			\sup_{x_1,x_2 \in B_{\eps^{-1}}(p)} \vert d(f(x_1),f(x_2)) - d(x_1,x_2) \vert \leq\epsilon, 
		\end{equation}
		\begin{equation}\label{pgh3}
			\sup_{y \in B_{\eps^{-1}}(q) } \inf_{x \in B_{\eps^{-1}}(p)}  d(f(x),y) \leq\eps. 
		\end{equation}
	\end{definition}
	
	\begin{definition}\label{pmgh}
		\rm Let $(X_i,p_i)$ be a sequence of pointed proper geodesic spaces.\ We say that it \textit{converges in the pointed Gromov--Hausdorff sense} (or pGH-sense) to a pointed proper geodesic space $(X,p)$ if there is a sequence of pointed $\epsilon_i$-GHAs $\phi_i : X_i \to X$ with $\epsilon_i \to 0 $ as $i\to\infty$. 
		
		If in addition to that, $(X_i,d_i,\mathfrak{m}_i)$, $(X,d,\mathfrak{m})$ are metric measure spaces, the maps $\phi_i$ are Borel measurable and $\int_X f \cdot d(\phi_i)_{\#}\mathfrak{m}_i \to \int_X f \cdot d\mathfrak{m}$, for all continuous $f: X \to \mathbb{R}$ with compact support, then we say $(X_i,d_i, \mathfrak{m}_i,p_i)$ converges to $(X,d,\mathfrak{m},p)$ in the \textit{pointed measured Gromov--Hausdorff sense} (or pmGH-sense).
	\end{definition}
	
	\begin{remark}
		\rm The maps $\phi_i:X_i\to X$ above are called Gromov-Hausdorff approximations and for any $x\in X$, there is a sequence $x_i\in X_i$ such that $\phi_i(x_i)\to x$.
	\end{remark}
	\begin{remark}
		If there is a sequence of groups $\Gamma_i$ acting on $X_i$ by (measure preserving) isometries with diam$(X_i/\Gamma_i) \leq D$ for some $D>0$, then one could ignore the points $p_i$ when one talks about p(m)GH convergence, as any pair of limits are going to be isomorphic as metric (measure) spaces.\ In particular, if $\op{diam}(X_i)\leq D$ for some $D>0$, then we simply say that $X_i$ converges to $X$ in the (m)GH-sense.
	\end{remark}
	
	One of the main features of the RCD($K,N$) condition is that it is stable under pmGH convergence, i.e., the pmGH-limit of a sequence of RCD($K,N$) spaces is also an RCD($K,N$) space (see \cite{gigli2015convergence}).\ Combined with Gromov's precompactness criterion and Prokhorov's compactness theorem (see \cite[Chapter 27]{villani2009optimal} for instance), the class of RCD($K,N$) spaces with normalized measure is compact under the pmGH-topology.
	
	\begin{definition}\label{def-tangent}
		We say that a pointed metric measure space $(Y,d_{Y},m_{Y},y)$ is a tangent cone of $(X,d,\mm)$ at $x$ if there exists a sequence $r_{i}\to0^{+}$ such that 
		\[(X,r_{i}^{-1}d,\mm(B_{r_{i}}(x))^{-1}\mm,x)\xrightarrow{pmGH}(Y,d_{Y},m_{Y},y).\]
		The collection of all tangent cones of $(X,d,\mm)$ at $x$ is denoted by $\text{Tan}_x(X,d,\mm)$.
	\end{definition}
	
	For an RCD$(K,N)$ space $(X,d,\mm)$, the compactness yields that $\text{Tan}_x(X,d,\mm)$ is non-empty and any tangent cone is RCD($0,N$).\ We are now in the position to introduce the notions of $k$-regular set and essential dimension as follows.
	
	\begin{definition}[$k$-regular set]\label{def-regular}
		For any integer $k\in[1,N]$, we denote by $\mathcal{R}_{k}$ the set of all points $x\in X$ such that $\text{Tan}_x(X,d,\mm)=\left\lbrace(\mathbb{R}^{k},d_{\mathbb{R}^{k}},(\omega_{k})^{-1}\mathcal{H}^{k},0^{k})\right\rbrace$, where $\omega_{k}$ is the volume of the unit ball in $\mathbb{R}^{k}$.\ We call $\mathcal{R}_{k}$ the $k$-regular set of $X$.
	\end{definition}
	\indent The following result is proved by Bru\`e-Semola in \cite{brue2020constancy}.
	\begin{theorem}[\citeonline{brue2020constancy}]\label{thm-essdim}
		Let $(X,d,\mm)$ be an RCD$(K,N)$ space with $K\in\mathbb{R}\ and\ N\in(1,\infty)$.\ Then there exists a unique integer $k\in[1,N]$, called the essential dimension of $(X,d,\mm)$, denoted by $\dim_{ess}(X)$, such that $\mm(X\setminus\mathcal{R}_{k})=0$.
	\end{theorem}

	\subsection{Equivariant Gromov--Hausdorff convergence}
	
	There is a well studied notion of convergence of group actions in this setting, first introduced by Fukaya-Yamaguchi in \cite{fukaya1992fundamental}.\ For a pointed proper metric space $(X,p)$, we equip its isometry group $\op{Iso}(X)$ with the metric $d_{0}^{p}$ given by
	\begin{equation}\label{d0}
		d_0^p (h_1, h_2) : = \inf_{r > 0 } \left\{ \frac{1}{r} + \sup_{x \in B_r(p)} d(h_1x, h_2x)  \right\}  .         
	\end{equation}
	for $h_1,h_2\in\op{Iso}(X)$.\ Obviously, we get a left invariant metric that induces the compact-open topology and makes $\op{Iso}(X)$ a proper metric space. 
	
	Recall that if a sequence of pointed proper metric spaces $(X_i, p_i)$ converges in the pGH sense to the pointed proper metric space $(X,p)$,  one has a sequence of pointed $\epsilon_i$-GHAs $\phi_i : X_i \to X$ with $\epsilon_i\to0$.
	
	\begin{definition}\label{def:equivariant}
		\rm Consider a sequence of pointed proper metric spaces $(X_i,p_i)$ that converges in the pGH sense to a pointed proper metric space $(X,p)$, a sequence of closed groups of isometries $G_i \leq \op{Iso}(X_i)$ and a closed group $G \leq \op{Iso}(X)$.\ Equip $G_i$ with the metric $d^{p_i}_{0}$ and $G$ with the metric $d^{p}_{0}$.\ We say that the sequence $G_i$ \textit{converges equivariantly} to $G$ if there is a sequence of Gromov-Hausdorff approximations $f_i:G_i\to G$ such that for each $R>0$, one has
		\[\lim\limits_{i\to\infty}\sup_{g \in B^{G_i}_R(Id_{X_i})}\sup_{x \in B^{X_i}_R(p_i)}d(\phi_i(gx),f_i(g)\phi_i(x))=0.\]
	\end{definition}

	Let us recall two basic properties of equivariant convergence proved in \cite{fukaya1992fundamental}.
	
	\begin{lemma}[\citeonline{fukaya1992fundamental}]\label{equivariant}
		Let $(Y_i,q_i)$ be a sequence of proper geodesic spaces that converges in the pGH sense to $(Y,q)$, and  $\Gamma_i \leq \op{Iso}(Y_i)$ a sequence of closed groups of isometries that converges  equivariantly to a closed group $\Gamma \leq \op{Iso}(Y)$.\ Then the sequence $(Y_i/\Gamma_i, [q_i])$ converges in the pGH sense to $(Y/\Gamma , [q])$.
	\end{lemma}

	\begin{lemma}[\citeonline{fukaya1992fundamental}]\label{equivariant-compactness}
		Let $(Y_i,q_i) $ be a sequence of proper geodesic spaces that converges in the pGH sense to $(Y,q)$, and take a sequence $\Gamma_i \leq \op{Iso}(Y_i)$ of closed groups of isometries.\ Then there is a subsequence $(Y_{i_k}, q_{i_k}, \Gamma_{i_k})_{k \in \mathbb{N}}$ such that $\Gamma_{i_k}$ converges equivariantly  to a closed group $\Gamma \leq \op{Iso}(Y)$.
	\end{lemma}
	A sequence of groups, which converges equivariantly to the trivial group, is called a sequence of small groups.\ The explicit definition is the following.
	\begin{definition}
		Let $(X_i, p_i)$ be a sequence of proper geodesic spaces. We say a sequence of groups $W_i \leq \iso (X_i) $ consists of \emph{small subgroups} if for each $R > 0 $ we have
		\[    \lim _{i \to \infty } \, \sup _{g \in W_i}\,  \sup _{x \in B_R(p_i)} d(gx,x) = 0 .         \]
		Equivalently, the groups $W_i$ are small if $d_0^{p_i}(g_i,Id_{X_i})\to0$ for any choice of $g_i\in W_i$.
	\end{definition}
	
	Obviously, if $W_i \leq \iso (X_i) $ is a sequence of discrete small subgroups, then $W_i$ is a finite group for any large $i$.\ It is proved in \cite[Theorem 93]{santos2023fundamental} that a non-collapsing sequence of $RCD(K,N)$ spaces cannot have small groups of measure preserving isometries.
	
	\begin{theorem}[\citeonline{santos2023fundamental}]\label{nss}
		Let $(X_i, d_i, \mathfrak{m}_i,p_i) $ be a sequence of pointed $RCD(K,N)$ spaces of essential dimension $n$ and let $H_i \leq Iso (X_i)$ be a sequence of small subgroups acting by measure preserving isometries.\ Assume the sequence $(X_i, d_i, \mathfrak{m}_i,p_i)$ converges in the pmGH sense to an $RCD(K,N)$ space $(X, d, \mathfrak{m},p) $ of essential dimension $n $. Then $H_i$ is trivial for large $i$.
	\end{theorem}

	\subsection{Properties of RCD(K,N) spaces} 
	
	In this subsection, we review some properties on RCD($K,N$) spaces that we will need in this paper. 
	
	Combining Theorem \ref{thm-essdim} and \cite[Theorem 4.1]{ambrosio2018short}, we have the following theorem.
	
	\begin{theorem}\label{thm-measure rectifiable}
		Let $(X,d,\mm)$ be an RCD$(K,N)$ space with $k=\dim_{ess}(X)$.\ Set 
		\[\mathcal{R}_k^*=\left\lbrace x\in\mathcal{R}_k :\ \lim\limits_{r\to0+}\frac{\mm(B_r(x))}{\omega_{k}r^k} \text{ exists and}\in(0,\infty)\right\rbrace. \]
		Then we have the following:
		\begin{enumerate}
			\item $\mm(X\setminus\mathcal{R}_k^*)=0$;
			\item $\mm\llcorner\mathcal{R}_k^*$ and $\mathcal{H}^k\llcorner\mathcal{R}_k^*$ are mutually absolutely continuous;
			\item $\lim\limits_{r\to0+}\dfrac{\mm(B_r(x))}{\omega_{k}r^k}=\dfrac{d\mm\llcorner\mathcal{R}_k^*}{d\mathcal{H}^k\llcorner\mathcal{R}_k^*}(x)$, for $\mm$-a.e.\ $x\in\mathcal{R}_k^*$.
		\end{enumerate}
	\end{theorem}
	
	The well known Cheeger--Gromoll splitting theorem \cite{cheeger1971splitting} was extended by Gigli to the setting  of RCD$(0,N)$ spaces \cite{gigli2014overview}.
	
	\begin{theorem}[\citeonline{gigli2014overview}]\label{gigli}
		Let $(X,d,\mathfrak{m})$ be an RCD$(0,N)$ space.\ Then there is a metric measure space $(Y, d_Y, \nu )$ where $(Y,d_Y)$ contains no line, such that $(X,d,\mathfrak{m})$ is isomorphic to the product $(\mathbb{R}^k\times Y, d_{\mathbb{R}^k} \times d_Y, \mathcal{H}^k \otimes \nu)$. Moreover, if $N-k \in [0,1)$ then $Y$ is a point, and in general, $(Y,d_Y,\nu)$ is an RCD$(0,N-k)$ space. 
	\end{theorem}
	
	\begin{remark}\label{rem-splitting}
		If $(X,d,\mm)$ is an RCD$(0,N)$ space and $G$ is a discrete subgroup of $\op{Iso}(X)$ with $\diam(X/G)<\infty$, then in the above splitting theorem, the space $Y$ can be taken to be compact (see \cite{mondino2019universal}).\ In addition, $\op{Iso}(X)=\op{Iso}(\mathbb{R}^k)\times \op{Iso}(Y)$.
	\end{remark}
	
	Let $(X,d,\mathfrak{m})$ be an RCD$(K,N)$ space  and $\rho : Y \to X$ be a covering space.\ $Y$ has a natural geodesic structure such that for any curve $\gamma : [0,1] \to Y$ one has 
	\[ \text{length}(\rho \circ \gamma) = \text{length}(\gamma).\]
	Set
	\[\mathcal{W} : = \{ W \subset Y  \text{ open bounded }\vert \text{ } \rho_{\vert_{W}} : W \to \rho (W)  \text{ is an isometry}  \}  \]
	and define a measure $\mathfrak{m}_{Y}$ on $Y$ by setting $\mathfrak{m}_Y(A) : = \mathfrak{m}(\rho (A))$ for each Borel set $A$ contained in $\mathcal{W}$.\ The measure $\mathfrak{m}_Y$ makes $\rho : Y \to X$ a local isomorphism of metric measure spaces, so by the local-to-global property of RCD$(K,N)$ space \cite{erbar2015equivalence}, $(Y,d_Y,\mathfrak{m}_Y)$ is an RCD$(K,N)$ space, and its group of deck transformations acts by measure-preserving isometries (see \cite{mondino2019universal} for more details). Wang proved the following in \cite{wang2024rcd}. 
	
	\begin{theorem}[\citeonline{wang2024rcd}]\label{jikang}
		Let $(X,d,\mathfrak{m})$ be an RCD$(K,N)$ space.\ Then for any $x\in X$ and $R>0$, there exists $r>0$ so that any loop in $B_r(x)$ is contractible in $B_R(x)$.\ In particular, $X$ is semi-locally-simply-connected and its universal cover $\tilde{X}$ is simply connected.
	\end{theorem}
	
	Due to Theorem \ref{jikang}, for an RCD$(K,N)$ space $(X,d,\mathfrak{m})$ we can think of its fundamental group $\pi_1(X)$ as the group of deck transformations of the universal cover $\tilde{X}$. 
	
	Recall that an RCD$(K,N)$ space $(X,d,\mm)$ is called non-collapsed if $\mm=\mathcal{H}^N$ (see \cite{de2018non}).\ There are some equivalent conditions for the non-collapse condition up to a scaling on measure (see \cite[Theorem 2.20]{brena2023weakly}, \cite[Theorem 2.3]{zamora2024topological} and references therein).
	\begin{theorem} \label{thm:dim}
		Let $(X,\mathsf{d},\mm)$ be an RCD$(K,N)$ space.\ Then the following five conditions are equivalent.
		\begin{enumerate}
			\item $X$ has essential dimension $N$. \label{thm:dim-1}
			\item $X$ has topological dimension $N$. \label{thm:dim-3}
			\item $X$ has Hausdorff dimension $N$. \label{thm:dim-2}
			\item $\mm = c \mathcal{H}^N$ for some constant $c>0$.\label{thm:dim-4}
			\item $N\in\mathbb{N}$ and $X$ has Hausdorff dimension greater than $N-1$.
		\end{enumerate}
	\end{theorem}
	
	De Philippis-Gigli \cite{de2018non} studied the GH-limit of non-collapsed RCD spaces and obtained the following result, generalizing Cheeger-Colding's result on Ricci-limit spaces \cite{cheeger1997structure}.
	
	\begin{theorem}[\citeonline{de2018non}]\label{thm:volume-continuity}
		Let $(X_i, d_i, \mathcal{H}^N, p_i )$ be a sequence of pointed  $\rcd (K,N)$ spaces and $(X_i,d_i,p_i)\xrightarrow{pGH}(X, d, p)$. Then precisely one of the following holds:
		\begin{enumerate}
			\item $\limsup_{i \to \infty}\mathcal{H}^N( B_1 (p_i) )>0$.\ In this case, $(X_i,d_i,\mathcal{H}^N,p_i)\xrightarrow{pmGH}(X, d,\mathcal{H}^N,p)$ and the limit $\lim_{i \to \infty}\mathcal{H}^N( B_1 (p_i) )$ exists and equal to $\mathcal{H}^N (B_1 (p))$. \label{eq:non-collapsed}
			\item $\lim_{i \to \infty}\mathcal{H}^N( B_1 (p_i) )=0$.\ In this case, $\dim_{\mathcal{H}}(X)\leq N-1$.
		\end{enumerate}
	\end{theorem}
	
	

	Let us recall the following topological stability theorem for non-collapsed RCD$(K,N)$ spaces, proved by Kapovitch-Mondino \cite[Theorem 3.3]{kapovitch2021topology}, based on Cheeger-Colding's Reifenberg type theorem \cite{cheeger1997structure}.
	
	\begin{theorem}[\citeonline{kapovitch2021topology}]\label{thm:reifenberg-weak}
		Let $(X_i,d_i,\mathcal{H}^N,p_i)$ be a sequence of pointed $\rcd (K,N)$ spaces such that the sequence $(X_i, p_i)$ converges in the pGH-sense to $(M^N,p)$ where $M^N$ is a smooth Riemannian manifold.\ Then for any $R>0$, there is a sequence of pointed $\epsilon_i$-GHAs $f_i : (X_i , p_i) \to ( M^N , p ) $ with $\epsilon_i\to0$, such that for all $i$ large enough depending on $R$, the restriction of $f_i$ to $B_R(p_i)$ is a bi-H\"older homeomorphism onto its image, and
		\[         B_{R - 4\epsilon_i}(p) \subset f_i (B_R(p_i)) .  \] 
		In particular, if $M$ is compact then $X_i$ is bi-H\"older homeomorphic to $M$ for all large $i$.
	\end{theorem}
	
	The bi-H\"older homeomorphism can be constructed via harmonic splitting map.\ We refer to \cite{brue2022boundary} for the notion of harmonic $(k,\eps)$-splitting map on RCD spaces.
	\begin{theorem}[\citeonline{cheeger2021rectifiability,honda2023note}]\label{Rei}
		Assume that $(X,d,\mathcal{H}^N,p)$ is a pointed RCD$(-\epsilon,N)$ space and $u: B_2(p) \to \mathbb{R}^N$ is a harmonic $(N,\epsilon)$-splitting map.\ Then for any $x,y \in B_1(p)$ we have 
		$$(1-\Phi(\epsilon|N))d(x,y)^{1+\Phi(\epsilon|N)} \le d(f(x),f(y)) \le (1+\Phi(\epsilon|N))d(x,y),$$
		where $\Phi(\eps|N)\to0$ as $\eps\to0$ while $N$ is fixed.\ Moreover, if $X$ is a smooth $N$-manifold with $\Ric\geq-\epsilon$, then for any $x \in B_1(p)$, $du:T_x X \to \mathbb{R}^N$ is nondegenerate.
	\end{theorem}

	\subsection{Nilpotent and polycyclic groups}\label{subsec:rank} 
	
	\begin{definition}
		For a group $G$, let $G^{(0)} : = G$ and define inductively $G^{(j + 1 )}:=[G^{(j)}, G]$.\ $G$ is called \emph{nilpotent} if $G^{(s)}$ is the trivial group for some $s \in \mathbb{N}$.\ $G$ is called \emph{virtually nilpotent} if there exists a nilpotent subgroup $N\leq G$ of finite index.
	\end{definition}
	\begin{definition}\label{def-polycyclic}
		A finitely generated group $\Lambda$ is said to be \emph{polycyclic} if there is a finite subnormal series
		\[  \Lambda= \Lambda _m \trianglerighteq  \cdots \trianglerighteq \Lambda_0 = {1}    \]
		with $\Lambda _ j / \Lambda _ {j-1} $ cyclic for each $j$.\ Such a subnormal series is called a polycyclic series.\ The polycyclic rank is defined as the number of $j$'s for which $\Lambda _ j / \Lambda _{j-1}$ is isomorphic to $\Z$, which is independent of the choice of the polycyclic series and denoted by $\rank (\Lambda )$.
	\end{definition}
	From the definition, we immediately know that any finite index subgroup of a polycyclic group is a polycyclic group with the same rank..
	
	It is well-known that any finitely generated nilpotent group is polycyclic (see \cite{kargapolov1979fundamentals} for instance), and the rank of a finitely generated nilpotent group is defined to be the polycyclic rank.\ The following lemma (see \cite[Lemma 2.22 and Lemma 2.24]{naber2016topology}) gives the definition of the rank of a finitely generated virtually nilpotent group.

	\begin{lemma}\label{lem:rank of virtually nilpotent}
		Let $G$ be a finitely generated virtually nilpotent group.\ Then:
		\begin{enumerate}
			\item Every nilpotent subgroup $N\leq G$ of finite index has the same rank.\ The common rank is called the rank of $G$ and also denoted by $\op{rank} (G)$.
			\item If $\Gamma$ is a finite index subgroup of $G$, then $\op{rank} (\Gamma)=\op{rank}(G)$.
		\end{enumerate}
	\end{lemma}
	
	Following \cite{zamora2024topological}, we can define the rank for an arbitrary finitely generated group.  
	\begin{definition}\label{def-rank}
		For a finitely generated group $G$, we define 
		\[\rank (G):=\inf\left\lbrace \rank (\Lambda ):\ \Lambda \text{ is a finite index polycyclic subgroup of } G \right\rbrace. \]
		The infimum of the empty set is defined to be $+\infty$.
	\end{definition}
	By Definition \ref{def-polycyclic} and Lemma \ref{lem:rank of virtually nilpotent}, if $G$ is polycyclic or finitely generated virtually nilpotent, there is no conflict between the distinct definitions of $\rank ( G )$.\ Also, if $\Lambda$ is a finite index subgroup of $G$, then $\op{rank}(\Lambda)=\op{rank}(G)$.
	
	We also note that if $G$ is a discrete group of isometries on a proper geodesic space $X$ with $\diam(X/G)\in(0,\infty)$, then by \cite[Lemma 2.5]{zamora2024limits}, $G$ is finitely generated.
	
		\subsection{Riemannian orbifolds}\label{subsec:orb}
	In this subsection, we review the basic theory of orbifolds.\ An orbifold $\orb$ is, roughly speaking, a topological space that is locally homeomorphic to a quotient of $\mathbb{R}^n$ by some finite group.\ We recall the definitions from \cite{kleiner2011geometrization} (see also \cite{galaz2018quotients}).
	
	
	\begin{definition}
		\label{D:local_model}
		A \emph{local model of dimension $n$} is a pair $(\hat{U}, G)$, where $\hat{U}$ is an open, connected subset of a Euclidean space $\RR^n$, and $G$ is a finite group acting smoothly and effectively on $\hat{U}$. 
		
		A \emph{smooth map} $(\hat{U}_1,G_1)\to (\hat{U}_2,G_2)$ between local models $(\hat{U}_i,G_i)$, $i=1,2$, is a homomorphism $f_{\sharp}: G_1\to G_2$ together with a $f_{\sharp}$-equivariant smooth map $\hat{f}:\hat{U}_1\to \hat{U}_2$, i.e., $\hat f(\gamma\cdot \hat u) = f_\sharp(\gamma) \cdot \hat f(\hat u)$, for all $\gamma \in G_1$ and $\hat u \in \hat U_1$.
	\end{definition}
	
	Given a local model $(\hat{U},G)$, denote by $U$ the quotient $\hat{U}/G$.\ The smooth map between local models is called an \emph{embedding} if $\hat{f}$ is an embedding. In this case, the effectiveness of the actions in the local models implies that $f_{\sharp}$ is injective.

	
	\begin{definition}
		An \emph{$n$-dimensional orbifold local chart} $(U_x, \hat U_x, G_x, \pi_x)$ around a point $x$ in a topological space $X$ consists of:
		\begin{enumerate}
			\item A neighborhood $U_x$ of $x$ in $X$;
			\item A local model $(\hat{U}_x, G_x)$ of dimension $n$;
			\item A $G_x$-equivariant projection $\pi_x:\hat{U}_x\to U_x$ that induces a homeomorphism $\hat{U}_x/G_x\to U_x$.
		\end{enumerate}
		If $\pi_x^{-1}(x)$ consists of a single point, $\hat{x}$, then $(U_x, \hat U_x, G_x, \pi_x)$ is called a \emph{good local chart} around $x$.\ In particular, $\hat x$ is fixed by the action of $G_x$ on $\hat{U}_x$.
	\end{definition}

	
	\begin{definition}
		An \emph{$n$-dimensional orbifold atlas} for a topological space $X$ is a collection of $n$-dimensional local charts $\mathcal{A}=\{U_{\alpha}\}_\alpha$ such that the local charts $U_\alpha \in \mathcal{A}$ give an open covering of $X$ and for any $x\in U_{\alpha}\cap U_{\beta}$, there is a local chart $U_\gamma \in \mathcal{A}$ with $x\in U_{\gamma}\subset U_{\alpha}\cap U_{\beta}$ and embeddings $(\hat{U}_{\gamma}, G_{\gamma})\to(\hat{U}_{\alpha}, G_{\alpha})$, $(\hat{U}_{\gamma}, G_{\gamma})\to (\hat{U}_{\beta}, G_{\beta})$.
		
		Two $n$-dimensional atlases are called \emph{equivalent} if they are contained in a third atlas. 
	\end{definition}

	
	\begin{definition}
		An \emph{$n$-dimensional (smooth) orbifold}, denoted by $\orb^n$ or simply $\orb$, is a second-countable, Hausdorff topological space $|\orb|$, called the \emph{underlying topological space} of $\orb$, together with an \emph{equivalence class of $n$-dimensional orbifold atlases}.  
	\end{definition}

	
	Given an orbifold $\orb$ and any point $x\in |\orb|$, one can always find a good local chart $U_x$ around $x$.  Moreover, the corresponding group $G_x$ does not depend on the choice of good local chart around $x$, and is referred to as the \emph{local group at $x$}.\ From now on, only good local charts will be considered.
	
	Each point $x\in |\orb|$ with $G_{x}=\{e\}$ is called a \emph{regular point}.\ The subset $|\orb|_{reg}$ of regular points is called \emph{regular part};  it is a a smooth manifold that forms an open dense subset of $|\orb|$.  A point which is not regular is called \emph{singular}.
	
	If a discrete group $\Gamma$ acts properly discontinuously on a
	manifold $M$, then the quotient space $M/\Gamma$ can be naturally endowed with an orbifold structure.\ For simplicity, we still use the terminology $M/\Gamma$ to mean $M/\Gamma$ as an orbifold.\ An orbifold $\Or$ is {\em good} if $\Or = M/\Gamma$ for some manifold $M$ and some discrete group $\Gamma$.
	
	Similarly, suppose that a discrete group $\Gamma$ acts 
	by diffeomorphisms on an orbifold 
	$\Or$.\ We say that it acts {\em properly discontinuously} if the action of $\Gamma$ on $|\Or|$ is properly discontinuous. Then there is a quotient orbifold $\Or/\Gamma$, with
	$|\Or/\Gamma| = |\Or|/\Gamma$.
	
	A {\em smooth map} $f  :  \Or_1 \rightarrow \Or_2$ between orbifolds
	is given by a
	continuous map $|f|  :  |\Or_1| \rightarrow |\Or_2|$ with the 
	property that for each $p \in |\Or_1|$, there are local models
	$(\hU_1, G_1)$ and $(\hU_2, G_2)$ for $p$ and $f(p)$ respectively, and a smooth map $\hat{f} : (\hU_1, G_1) \rightarrow (\hU_2, G_2)$ between local models so that the diagram
	\begin{equation}
		\begin{matrix}
			\hU_1 & \stackrel{\hat{f}}{\rightarrow} & \hU_2 \\
			\downarrow & & \downarrow \\
			U_1 & \stackrel{|f|}{\rightarrow} & U_2
		\end{matrix}
	\end{equation}
	commutes.\ A {\em diffeomorphism} $f \: : \: \Or_1 \rightarrow \Or_2$ is a smooth map with a smooth inverse.\ In this case, $G_p$ is isomorphic to $G_{f(p)}$.
	

	An {\em orbifold covering} $\pi : \hat\Or \rightarrow \Or$ is a surjective smooth map such that 
	\begin{enumerate}
		\item for each $x\in|\orb|$, there is an orbifold local chart $(U,\tilde U,H,\phi)$ around $x$ such that $|\pi|^{-1}(U)$ is a disjoint union of open subsets $V_i\subset |\hat{\Or}|$;
		\item each $V_i$ admits an orbifold local chart of the type $(V_i,\tilde U,H_i,\phi_i)$ where $H_i<H$, such that $|\pi|$ locally lifts to the identity $\tilde{U}\to\tilde{U}$ with inclusion $H_i\to H$.
	\end{enumerate}
	A universal orbifold covering of $\Or$ is an orbifold covering $\pi : \tilde\Or \rightarrow \Or$ such that for every orbifold covering $\psi: \hat\Or \rightarrow \Or$, there is an orbifold covering $\phi: \tilde\Or \rightarrow \hat\Or$ so that $\psi\circ\phi=\pi$.\ It is due to Thurston \cite{thurston1997three} that any connected orbifold $\Or$ admits a universal orbifold covering $\pi : \hat\Or \rightarrow \Or$.\ The {\em orbifold fundamental group} of $\Or$, denoted by $\pi_1^{orb}(\Or)$, is defined to be the deck transformation group of its universal orbifold covering.\ The universal orbifold covering $\pi : \tilde\Or \rightarrow \Or$ induces a diffeomorphism $\tilde\Or /\pi_1^{orb}(\Or) \to \Or$.
	
	In general, an orbifold covering is not locally homeomorphism and hence, not a covering.\ In addition, $\pi_1^{orb}(\Or)$ is different from $\pi_1(|\Or|)$ and there is actually a epimorphism $\pi_1^{orb}(\Or)\to \pi_1(|\Or|)$ (see \cite{haefliger1984appendice}).
	
	

	\begin{definition} [Riemannian metric on an orbifold]
		A Riemannian metric $g$ on an orbifold $\orb$ is given by a collection of Riemannian metrics on the local models $\hat{U}_\alpha$ so that the following conditions hold:
		\begin{enumerate}
			\item The local group $G_\alpha$ acts isometrically on $\hat{U}_\alpha$.
			\item The embeddings $(\hat{U}_3,G_3)\to (\hat{U}_1,G_1)$ and $(\hat{U}_3,G_3)\to (\hat{U}_2,G_2)$ in the definition of orbifold atlas are isometries (with respect to the Riemannian metric).
		\end{enumerate}
	\end{definition}
	
	Note that the Riemannian metric $g$ induces a natural metric $d$ on $|\orb|$ that is locally isometric to the quotient metric of $(\hat{U}_{x},\hat{d}_{x})$ by $G_{x}$, where $\hat{d}_{x}$ is induced by the Riemannian metric $\hat g_x$ on $\hat{U}_{x}$.\ In the absence of ambiguity, we sometimes directly treat $(\orb,g)$ as the metric space $(|\orb|,d)$ and apply the terminology from metric spaces to $(\orb,g)$.
	
	
	For any Riemannian  orbifold $(\orb,g)$, there is a natural volume  measure $\operatorname{vol}_{g}$ given on the local orbifold charts $(U_x, \hat U_x, G_x, \pi_x)$ by 
	$\operatorname{vol}_{g}|_{U_{x}}:= \frac{1}{|G_x|}(\pi_{x})_{\#}\operatorname{vol}_{\hat g_{x}}$,
	where $\operatorname{vol}_{\hat g_{x}}$ is the Riemannian volume measure on $(\hat{U}_{x}, \hat g_{x})$.
	
	The regular part $|\Or|_{reg}$ inherits a Riemannian metric.\ The corresponding volume form equals the $n$-dimensional Hausdorff measure on $|\Or|_{reg}$. In particular, $\op{vol}_g(\Or)$ coincides with the volume of the Riemannian manifold $|\Or|_{reg}$, which equals the $n$-dimensional Hausdorff measure of the metric space $|\Or|$.
	
	
	
	The Levi-Civita connection on $(\orb,g)$ can be defined via the local models.\ We can then define the curvature tensor $Rm$ on $\orb$ and derived curvature notions, such as sectional and Ricci curvatures, are defined accordingly.\ We adopt the same notation for corresponding geometric quantities as is used on Riemannian manifolds.
	
	Letting $g_{reg} = g|_{|\orb|_{reg}}$, we have that 
	$(|\orb|_{reg},g_{reg})$ is a smooth open Riemannian manifold.\ By density, it is clear that $(|\orb|_{reg},g_{reg})$ satisfies $\Ric_{g_{reg}}\geq K$ if and only if $(\orb, g)$ satisfies $\Ric_{g}\geq K$.
	
	Also, the following result was proved by Galaz-Kell-Mondino-Sosa \cite[Theorem 7.10]{galaz2018quotients}.
	
	\begin{theorem}[\citeonline{galaz2018quotients}]\label{thm-orb rcd}
		Let $(\orb,g)$ be an $n$-dimensional Riemannian orbifold.\ Then $\Ric_g\geq K$ if and only if $(\orb,g,\op{vol}_g)$ is an RCD$(K,n)$ space.
	\end{theorem}

	Finally, let us review some facts about closed flat orbifolds.\ Recall that a group $\Gamma\leq \op{Iso}(\R^n)$ is called \emph{crystallographic} if it is discrete and cocompact, so that $\R^n/\Gamma$ is a closed flat orbifold.\ Conversely, by a result of Thurston \cite{thurston1997three}, if $(\orb,g)$ is a closed flat orbifold, then it is good, its universal orbifold cover is $\R^n$ and $\pi_1^{orb}(\orb)$ is isomorphic to a crystallographic group.

	
	
	\subsection{Almost-crystallographic groups and infranil orbifolds}
	
	In this subsection, we review some basic notions of almost-crystallographic groups and infranil orbifolds.
	
	\begin{definition}\label{def-almost crystal}
		Let $\mathcal{N}$ be a connected and simply connected nilpotent Lie group, and consider a maximal compact subgroup $C$ of $\op{Aut}(\mathcal{N})$.\ A cocompact and discrete subgroup $\Gamma$ of $\mathcal{N}\rtimes\op{Aut}(\mathcal{N})$ is called an almost-crystallographic group (modeled on $\mathcal{N}$).\ The dimension of $\Gamma$ is defined to be that of $\nn$.\ If moreover, $\Gamma$ is torsion free, then $\Gamma$ is called an almost-Bieberbach group.
		
		An infranil orbifold (resp.\ infranilmanifold) is a quotient space $\mathcal{N}/\Gamma$, where $\Gamma$ is an almost-crystallographic (resp.\ almost-Bieberbach) group modeled on $\mathcal{N}$.\ If further $\Gamma\subset\mathcal{N}$ (so $\Gamma$ acts freely on $\mathcal{N}$), then we say that $\mathcal{N}/\Gamma$ is a nilmanifold.
	\end{definition}
	
	Let $\Gamma$ be an almost-crystallographic group modeled on $\mathcal{N}$.\ Due to the generalized first Bieberbach Theorem proved by Auslander \cite{auslander1960bieberbach}, $G=\nn\cap\Gamma$ is a lattice of $\nn$ and $\Gamma/G$ is finite.\ Therefore, an infranil orbifold is actually the quotient of a nilmanifold by a finite group of affine diffeomorphisms.\ In addition, $\op{rank}(\Gamma)=\op{rank}(G)=\dim(\nn)$.
	
	Let us recall a well-known algebraic characterization of almost-crystallographic groups (see \cite[Theorem 4.2]{dekimpe2016users} for instance).
	
	\begin{theorem}\label{thm-almost crystallographic}
		Let $E$ be a finitely generated virtually nilpotent group.\ Then the followings are equivalent.
		\begin{enumerate}
			\item $E$ is isomorphic to an almost-crystallographic group.
			\item $E$ contains a torsion free nilpotent normal subgroup $G$, such that $G$ is maximal nilpotent in $E$ and $[E:G]<\infty$.
			\item $E$ does not contain a non-trivial finite normal subgroup.
		\end{enumerate}
	\end{theorem}
	
	To prove “(3) implies (1)" part in the above theorem, one may first find a normal subgroup $N$ of finite index in $E$ such that $N$ is a finitely generated torsion free nilpotent group.\ It was shown by Mal'cev \cite{malcev1949} that such $N$ can always be embedded as a lattice in a simply connected nilpotent Lie group $\nn$, which is unique up to isomorphism and now called the Mal'cev completion of $N$.\ Then $E$ can be identified to an almost-crystallographic group modeled on $\nn$.\ Specifically, we have the following proposition (see \cite{dekimpe2016users} for the detailed proof).
	
	\begin{proposition}\label{prop-Malcev}
		Let $E$ be a finitely generated virtually nilpotent group which does not contain a non-trivial finite normal subgroup.\ Let $N$ be a torsion free nilpotent normal subgroup of finite index in $E$.\ Then $E$ is isomorphic to an almost-crystallographic group modeled on $\nn$, where $\nn$ is the Mal'cev completion of $N$.
	\end{proposition}

	\subsection{Limits of almost homogeneous spaces}
	Let $(X_i,p_i)$ be a sequence of pointed almost homogeneous spaces, which converges in the pGH sense to $(X,p)$.\ By Lemma \ref{equivariant} and Lemma \ref{equivariant-compactness}, there is a closed group $G\leq \iso(X)$ acting transitively on $X$; that is, $X$ is $G$-homogeneous.\ Indeed, the limit of almost homogeneous spaces was specifically studied by Zamora in \cite{zamora2024limits}, where he utilized the results of Breuillard-Green-Tao \cite{breuillard2012structure} and proved the following theorem.
	
	\begin{theorem}[\citeonline{zamora2024limits}]\label{thm-zamora}
		Let $(X_i,p_i)$ be a sequence of pointed almost homogeneous spaces, which converges in the pGH sense to $(X,p)$.\ If $X$ is semi-locally-simply-connected, then $X$ is a nilpotent Lie group equipped with a sub-Finsler invariant metric, and $\pi_1(X)$ is a torsion free subgroup of a quotient of $\pi_1(X_i)$ for sufficiently large $i$.
	\end{theorem}
	Indeed, the fundamental group of any connected nilpotent Lie group is finitely generated torsion free abelian (see \cite[Corollary 2.11]{zamora2024limits}).\ In the above theorem, $X$ will be simply connected when $\pi_1(X_i)$ are finite groups.
	
	It is well-known that any compact connected nilpotent Lie group is abelian (see \cite[Corollary 2.13]{zamora2024limits} for instance).\ Thus, if the above limit space $X$ is compact, then it must be a torus.\ We note that this result on compact limits of almost homogeneous spaces can also be obtained in the finite dimensional case by using an old theorem of Turing \cite{turing1938finite}, and Gelander \cite{gelander2013limits} proved a more general result which covers the infinite dimensional case (see also \cite[Theorem 1.4]{zamora2024limits}).
	
	
	A key to proving Theorem \ref{thm-zamora} lies in finding a nilpotent group of isometries acting transitively on $X$, which was further applied in the recent work of Zamora-Zhu \cite{zamora2024topological}.\ The following three lemmas stem from \cite{zamora2024topological}.
	
	\begin{lemma}[\citeonline{zamora2024topological}]\label{lem:almost-homogeoeus-subgroups}
		Let  $(X_i, p_i)$ be a sequence of pointed proper geodesic spaces that converges to a pointed proper semi-locally-simply-connected geodesic space $(X,p)$ in the pointed Gromov--Hausdorff sense, and $G_i \leq \iso (X_i)$ a sequence of discrete groups with $\diam (X_i / G_i ) \to 0$.\ Then there exists $s \in \mathbb{N}$ and a sequence of finite index normal subgroups $ G_i  ^{\prime} \leq G_i $ with 
		$$\lim\limits_{i \to \infty }  \, \sup _{x \in X_i} \, \sup _{g \in (G_i ^{\prime } ) ^{(s)}}  d(gx, x )  = 0  \text{ and } \limsup\limits_{i \to \infty} \, [G_i : G_i ^{\prime}]  < \infty. $$ 
	\end{lemma}
	
	Note the basic fact that any subgroup of bounded index contains a normal subgroup of bounded index (see \cite[Lemma 4.8]{wang2024non} for instance).\ Thus based on \cite[Lemma 2.23]{zamora2024topological}, we can further assume that $G_i^{\prime}$ is a normal subgroup.
	
	\begin{lemma}[\citeonline{zamora2024topological}]\label{lem:almost-translations}
		Let $(X_i, p_i)$, $(X,p)$, $G_i$, $G_i^{\prime}$ be as in Lemma \ref{lem:almost-homogeoeus-subgroups}.\ Then after passing to a subsequence the groups $G_i^{\prime}$ converge equivariantly to a connected nilpotent group $ G\leq \iso (X)$ acting freely and transitively.  
	\end{lemma}
	
	\begin{lemma}[\citeonline{zamora2024topological}]\label{lem:free}
		Let $(X_i, p_i)$, $(X,p)$, $G_i$, $G_i^{\prime}$ be as in Lemma \ref{lem:almost-homogeoeus-subgroups}.\ If $G_i$ satisfy that any sequence of small subgroups is trivial for large $i$, then $G_i^{\prime}$ acts freely for large $i$. 
	\end{lemma}
	
	We further introduce the following definition and lemma from \cite{breuillard2012structure,zamora2024limits,zamora2024topological}, which gives an explicit description of the groups $G_i^{\prime}$.

	\begin{definition}
		Let $G$ be a group, $u_1, u_2, \ldots , u_r \in G$, and $N_1, N_2, \ldots , N_r $ $\in \mathbb{R}^+$. The set $P(u_1, \ldots , u_r ; N_1, \ldots , N_r) \subset G$ is defined to be the set of elements that can be expressed as words in the $u_i$'s and their inverses such that the number of appearances of $u_i$ and $u_i^{-1}$ is not more than $N_i$.  We then say that  $P(u_1, \ldots , u_r ; N_1, \ldots , N_r)$  is a \emph{nilprogression in $C$-normal form} for some $C>0$ if it also satisfies the following properties:
		\begin{enumerate}
			\item  For all $1 \leq i \leq j \leq r$, and all choices of signs, we have
			\begin{center}
				$  [ u_i^{\pm 1} , u_j^{\pm 1} ]  \in P \left(  u_{j+1} , \ldots , u_r ; \dfrac{CN_{j+1}}{N_iN_j} , \ldots , \dfrac{CN_r}{N_iN_j}  \right). $
			\end{center}\label{condition:n1}
			\item The expressions $ u_1 ^{n_1} \ldots u_r^{n_r} $ represent distinct  elements as $n_1, \ldots , n_r$ range over the integers with  $\vert n_1 \vert \leq  N_1/C , \ldots , \vert n_r \vert \leq  N_r/C$.\label{condition:n2}
			\item One has 
			\[  \frac{1}{C} (2\lfloor N_1 \rfloor + 1 ) \cdots ( 2\lfloor N_r \rfloor + 1 ) \leq \vert P \vert  \leq C (2 \lfloor N_1 \rfloor + 1 ) \cdots ( 2\lfloor N_r \rfloor + 1 ).\] \label{condition:n3}
		\end{enumerate} 
		For a nilprogression $P$ in $C$-normal form, and $\varepsilon \in (0,1)$, the set  $P( u_1, \ldots , u_r ; $ $ \varepsilon N_1, \ldots , \varepsilon N_r   )$ also satisfies conditions (1) and (2), and we denote it by $\varepsilon P$.\ We define the \emph{thickness} of $P$ as the minimum of $N_1, \ldots , N_r$ and we denote it by $ \thi (P)$.\ The set $\{ u_1 ^{n_1 }\ldots u_r^{n_r} \vert \vert n_i \vert \leq N_i/C  \}$ is called the \emph{grid part of }$P$, and is denoted by $G(P)$.
	\end{definition}

	\begin{lemma}\label{lem:malcev-construction}
		Let $(X_i, p_i)$, $(X,p)$, $G_i$, $G_i^{\prime}$ be as in Lemma \ref{lem:almost-homogeoeus-subgroups} and $n=\dim_{top}(X)$.\ Then for $i$ large enough, there are $N_{1,i}, \ldots , N_{n,i} \in \mathbb{R}^+$ and torsion free nilpotent groups $\tilde{\Gamma}_i$ generated by elements $\tilde{u}_{1,i}, \ldots , \tilde{u}_{n,i} \in \tilde{\Gamma}_i$ with the following properties:
		\begin{enumerate}
			\item There are polynomials $Q_i : \mathbb {R}^n \times \mathbb{R}^n\to \mathbb{R}^n$ of degree $\leq d(n)$ giving the group structures on $\R^n$ by $x_1 \cdot x_2 = Q_i (x_1, x_2)$ such that for each $i$, $\tilde{\Gamma}_i$ is isomorphic to the group $(\Z^n,Q_i|_{\Z^n\times\Z^n})$ and the group $(\R^n,Q_i)$ is isomorphic to the Mal'cev completion of $\tilde{\Gamma}_i$.
			\label{conclusion:malcev-1}
			\item  There are small normal subgroups $W_i \trianglelefteq G_i^{\prime }$ and surjective group morphisms 
			\[\Phi _i : \tilde {\Gamma }_ i \to \Gamma_i : = G_i^{\prime } / W_i\]
			such that $\op{Ker} (\Phi_i)$ is a quotient of $\pi_1(X_i)$ and contains an isomorphic copy of $\pi_1(X)$ for $i$ large enough.  \label{conclusion:malcev-2}
			\item There is $C > 0 $ such that if $u_{j,i} : = \Phi _i (\tilde{u}_{j,i}) $ for each $j \in \{ 1, \ldots , n\}$, the set
			\[   P_i : = P (  u_{1,i}, \ldots , u_{n,i} ; N_{1,i}, \ldots , N_{n,i} )  \subset \Gamma_i       \]
			is a nilprogression in $C$-normal form with $\thi (P_i ) \to \infty $.  \label{conclusion:malcev-3}
			\item  \label{conclusion:malcev-4} For each $\varepsilon > 0 $ there is $\delta > 0 $ such that 
			\begin{align*}
				G(\delta P_i )     \subset \{ g \in \Gamma_i \, \vert  \,  &d(g[p_i], [p_i] ) \leq \varepsilon \}  , \\
				\{ g \in \Gamma_i \, \vert  \, d(g [p_i ], [p_i&] ) \leq \delta \}  \subset     G( \varepsilon  P_i   ) 
			\end{align*}
			for $i$ large enough, where we are considering the action of $\Gamma _ i$ on $X_i / W_i $. 
		\end{enumerate}
	\end{lemma}
	
	\begin{proof}
		 Note that \cite[Lemma 2.30]{zamora2024topological} provides (2), (3), (4), except that in (2) we additionally obtain that $\op{Ker} (\Phi_i)$ is a quotient of $\pi_1(X_i)$.\ This is due to \cite[Proposition 8.4]{zamora2024limits}.\ Moreover, (1) is just the Mal'cev Embedding Theorem and the construction of the polynomial $Q_i$ can be found in \cite[Section 5.1]{buser1981gromov} (see also \cite[Section 8]{zamora2024limits}).
	\end{proof}

	\begin{remark}\label{rem:rank-inequality}
		As noted in \cite[Remark 2.31]{zamora2024topological}, one obtains that
		\begin{equation}\label{eq:rank-string}
				\rank (G_i^{\prime} )  =  \rank (\Gamma _i ) = \rank (\tilde {\Gamma}_i) - \rank (\Ker (\Phi _ i)) \leq  n.
		\end{equation}
		If $\op{rank}(G_i^{\prime})=n$, then $\Ker (\Phi _ i ) $ is trivial and hence, $X$ is simply connected and $\tilde {\Gamma }_i = \Gamma_i$.
	\end{remark}
	
	\begin{remark}\label{rem-simply connected}
		If $\pi_1(X_i)$ are finite groups, then $\Ker(\Phi_i)$ will also be trivial, so $X$ is simply connected and $\tilde {\Gamma }_i = \Gamma_i= G_i^{\prime } / W_i$.\ In addition, $\rank(G_i)=\rank (G_i^{\prime})=\rank(\tilde {\Gamma}_i)=n$.
	\end{remark}
	
	\subsection{Topological rigidity for RCD spaces with bounded covering geometry}
	
	As we have noted in Section \ref{section-0}, Theorem \ref{thm-almostflat} can be extended to the RCD setting as the following theorem by the recent works of Zamora-Zhu \cite{zamora2024topological} and Wang \cite{wang2024non}.
	
	\begin{theorem}[Theorem \ref{thm-rcd-almostflat}]\label{yyc}
		For each $K\in\mathbb{R}$ and $N\geq1$, there is $\epsilon=\epsilon(K,N), v=v(K,N)$ such that for any RCD($K,N$) space $(X,d,\mathcal{H}^N)$ with $\op{diam}(X)\leq\epsilon$, the followings are equivalent:
		\begin{enumerate}
			\item $X$ is bi-H\"older homeomorphic to an $N$-dimensional infranilmanifold $\mathcal{N}/\Gamma$, where $\mathcal{N}$ is a simply connected nilpotent Lie group endowed with a left invariant metric;
			\item $\op{rank}(\pi_1(X))$ is equal to $N$;
			\item $\mathcal{H}^N(B_1(\tilde x))\geq v$, where $\tilde x$ is a point in the universal cover $\widetilde{X}$.
		\end{enumerate}
	\end{theorem}
	
	In the above theorem, (1) trivially implies (2).\ It was proved by Zamora-Zhu \cite{zamora2024topological} that (2) implies $X$ being homeomorphic to an $N$-dimensional infranilmanifold and their proof implicitly shows that (2) implies (3).\ In \cite{wang2024non}, Wang proved that (3) implies (1).
	
    We say that a non-collapsed RCD$(K,N)$ space $(X,d,\mathcal{H}^N)$ satisfies (global) $(1,v)$-bounded covering geometry if condition (3) of Theorem \ref{yyc} holds.
	
	Indeed, a local version of this term was proposed by Rong on Riemannian manifolds.\ Specifically, let $M$ be a compact $n$-manifold with $\Ric_M\geq-(n-1)$.\ We say that $M$ satisfies local $(r,v)$-bounded covering geometry if for any $x\in M$, $\op{vol}(B_{r/2}(\tilde{x}))\geq v$ where $\tilde{x}$ is a pre-image of $x$ in the (incomplete) Riemannian universal covering 
	\[\pi: (\widetilde{B_r(x)},\tilde{x})\to (B_r(x),x).\]
	If $M$ has small diameter, then local bounded covering geometry is equivalent to (global) bounded covering geometry.\ We refer to \cite{huang2020collapsed,rong2022collapsed} and the survey paper \cite{huang2020collapsing} for more detailed descriptions on this terminology.
	
	Naturally, one can define a non-collapsed RCD$(K,N)$ space $(X,d,\mathcal{H}^N)$ satisfying local bounded covering geometry in a similar fashion.\ Note that by Theorem \ref{jikang}, any $r$-ball $B_r(x)\subset X$ is semi-locally-simply-connected and hence, admits a universal cover. 
	
	The following fibration theorem summarizes the contributions from \cite{huang2020fibrations} and \cite{rong2022collapsed}.
	\begin{theorem}[\citeonline{huang2020fibrations,rong2022collapsed}]\label{thm-local covering geometry}
		Let $(M_i^n,g_i)$ converge to $(N^k,g)$ in the pGH-sense, where $N^k$ is a compact smooth manifold with $k\leq n$.\ Suppose that $(M_i^n,g_i)$ satisfy $\Ric_{g_i}\geq-(n-1)$ and local $(1,v)$-bounded covering geometry for some $v>0$.\ Then for all large $i$, there exist smooth fiber bundle maps $f_i: M_i\to N$ which are also $\epsilon_i$-GHAs with $\epsilon_i\to0$.\ Moreover, any $f_i$-fiber is diffeomorphic to an $(n-k)$-dimensional infranilmanifold.
	\end{theorem}
	
	In \cite{huang2020fibrations}, Huang constructed the smooth fiber bundle map and Rong \cite{rong2022collapsed} further identified the fibers to infranilmanifolds.\ This theorem is a generalization of Fukaya's fibration theorem in \cite{fukaya1987collapsing} on collapsed manifolds with bounded sectional curvature (see also \cite{cheeger1992nilpotent}).\ Indeed, any $n$-manifold with $|\sec|\leq1$ satisfies local $(r,v)$-bounded covering geometry with $r$ and $v$ depending only on $n$ (see \cite{cheeger1992nilpotent}).
	
	
	Recently, Wang \cite{wang2024non} generalized Theorem \ref{thm-local covering geometry} to RCD$(K,N)$ spaces $(X,d,\mathcal{H}^N)$.
	
	\begin{theorem}[\citeonline{wang2024non}]\label{wang-fibration}
		Let $(X_i,d_i,\mathcal{H}^N)$ be a sequence of RCD$(K,N)$ spaces and $(X_i,d_i)$ converge in the GH-sense to a closed Riemannian manifold $N^k$ with $k\leq N$.\ Suppose that all $(X_i,d_i,\mathcal{H}^N)$ satisfy local $(1,v)$-bounded covering geometry for some $v>0$.\ Then for large enough $i$, there are fiber bundle maps $f_i: X_i\to N^k$ which are also $\eps_i$-GHAs with $\eps_i\to0$ and any $f_i$-fiber with the induced metric is bi-H\"older homeomorphic to an $(N-k)$-dimensional infranilmanifold. 
	\end{theorem}

	\section{Limits of almost homogeneous RCD spaces and applications}\label{sec-3}
	
	In this section, we will prove Theorem \ref{thm-almost homogeneous rcd} and Theorem \ref{thm-almost nonnegative ricci-limit} and derive a series of consequences.\ The following lemma is obvious but key to our proof of Theorem \ref{thm-almost homogeneous rcd} (2).
	
	\begin{lemma}\label{lem-measure converge}
		Let $(X_i,d_i,m_i,q_i)$ be a sequence of metric measure spaces that converges in the pmGH sense to $(X,d,m,q)$, and $G_i \leq \op{Iso}(X_i)$ a sequence of closed groups of measure preserving isometries that converges equivariantly to a closed group $G \leq \op{Iso}(X)$.\ Then $G$ acts on $X$ by measure preserving isometries.
	\end{lemma}
	\begin{proof}
		Let $g$ be an arbitrary element of the group $G$.\ We need to show that $g_{\#}m=m$.
		
		Notice that by Definition \ref{def:equivariant}, there is a sequence of Gromov-Hausdorff approximations $f_i: G_i\to G$ and $g_i\in G_i$, such that $f_i(g_i)\to g$.\ Thus, $f_i(g_i)_{\#}m\to g_{\#}m$ in $C_{c}(X)^*$.\ Also, there is a sequence of Gromov-Hausdorff approximations $\phi_i: X_i\to X$ such that $(\phi_i)_{\#}m_i\to m$ in $C_{c}(X)^*$.\ Then we have $(f_i(g_i)\circ\phi_i)_{\#}m_i\to g_{\#}m$ in $C_{c}(X)^*$.
		
		By Definition \ref{def:equivariant}, $(\phi_i\circ g_i)_{\#}m_i\to g_{\#}m$ in $C_{c}(X)^*$.\ Since $(g_i)_{\#}m_i=m_i$, we obtain that $(\phi_i)_{\#}m_i\to g_{\#}m$ in $C_{c}(X)^*$.\ This leads to $g_{\#}m=m$ and we complete the proof.
	\end{proof}
	
	We say that a metric measure space $(X,d,m)$ is metric measure homogeneous if for all $x,y\in X$, there exists a measure preserving isometry $h\in \op{Iso}(X)$ such that $h(x)=y$.
	
	\begin{lemma}\label{lem-mm homogeneous}
		Let $(X,d,m)$ be a metric measure homogeneous RCD$(K,N)$ space of essential dimension $n$ for some $K\in\mathbb{R}$ and $N\in [1,\infty)$.\ Then $X$ is isometric to a Riemannian $n$-manifold and $m=c\mathcal{H}^n$ for some $c>0$.\ In particular, $(X,d,\mathcal{H}^n)$ is a non-collapsed RCD$(K,n)$ space.
	\end{lemma}
	\begin{proof}
		By \cite[Proposition 5.14]{santos2020invariant}, $X$ is isometric to a Riemannian $n$-manifold.\ Due to homogeneity, $X=\mathcal{R}_n^*$ (see Theorem \ref{thm-measure rectifiable}) and the limit
		\[\lim\limits_{r\to0^+}\frac{m(B_r(x))}{\omega_{n}r^n}\] 
		is a constant denoted by $c$.\ Then by Theorem \ref{thm-measure rectifiable}, we have $m=c\mathcal{H}^n$ and hence, $(X,d,\mathcal{H}^n)$ is RCD$(K,N)$.\ Since $X$ is a Riemannian $n$-manifold, $(X,d,\mathcal{H}^n)$ is an RCD$(K,n)$ space.\ This completes the proof.
	\end{proof}
	\begin{remark}
		In a recent work \cite{honda2024locally}, Honda-Nepechiy arrived at the same conclusion as Lemma \ref{lem-mm homogeneous} by assuming only that $(X,d,m)$ is \textit{locally} metric measure homogeneous.\ For our purposes, Lemma \ref{lem-mm homogeneous} is sufficient, and the proof is considerably simpler.
	\end{remark}
	
	Now, we can prove Theorem \ref{thm-almost homogeneous rcd} and for the convenience of readers, we rewrite it here.
	\begin{theorem}[Theorem \ref{thm-almost homogeneous rcd}]
		Let $K\in\mathbb{R}, N\in(1,\infty)$ and $(X_i,d_i,m_i,p_i)$ be a sequence of $(\epsilon_i,G_i)$-homogeneous pointed RCD$(K,N)$ spaces converging to $(X,d,m,p)$ in the pmGH-sense with $\epsilon_i\to0$.\ Assume that $X$ is not a point.\ Then the followings hold:
		\begin{enumerate}
			\item $(X,d)$ is isometric to an $n$-dimensional nilpotent Lie group with a left invariant Riemannian metric for some $n\leq N$;
			\item if the actions of $G_i$ are measure-preserving, then $m=c\mathcal{H}^n$ for some $c>0$ and $(X,d,\mathcal{H}^n)$ is an RCD$(K,n)$ space.\ In particular, $\op{Ric}_X\geq K$ if $n\geq2$;
			\item if $X_i$ are compact (or equivalently, $G_i$ are finite), then $m$ and $\mathcal{H}^n$ are mutually abolutely continuous and $(X,d,\mathcal{H}^n)$ is an RCD$(K,n)$ space.\ In particular, $\op{Ric}_X\geq K$ if $n\geq2$;
			\item if $X$ is compact, then $X$ is isometric to a flat torus $\mathbb{T}^n$.
		\end{enumerate}
	\end{theorem} 
	\begin{proof}
		(1) By Lemma \ref{equivariant} and Lemma \ref{equivariant-compactness}, $(X,d,m)$ is an homogeneous RCD$(K,N)$ space.\ Then by \cite[Proposition 5.14]{santos2020invariant}, $X$ is a Riemannian $n$-manifold, where $n=\dim_{ess}(X)\leq N$.\ Combining with Theorem \ref{thm-zamora}, $(X,d)$ is isometric to an $n$-dimensional nilpotent Lie group with a left invariant Riemannian metric.
		
		(2) By Lemma \ref{equivariant-compactness} and Lemma \ref{lem-measure converge}, there is a closed group $G\leq\op{Iso}(X)$ acting transitively on $X$ by measure preserving isometries.\ Then by Lemma \ref{lem-mm homogeneous}, $m=c\mathcal{H}^n$ for some $c>0$ and $(X,d,\mathcal{H}^n)$ is an RCD$(K,n)$ space.\ Since $X$ is a Riemannian n-manifold, $\op{Ric}_X\geq K$ if $n\geq2$.
		
		(3) Since $G_i$ is a finite group, we can apply \cite[Theorem A]{santos2020invariant} to obtain a $G_i$-invariant measure $m_{G_i}$ (so $G_i$ acts on $(X_i,d_i,m_{G_i})$ by measure preserving isometries) such that $(X_i,d_i,m_{G_i})$ is also RCD$(K,N)$.\ After normalizing the measure $m_{G_i}$ and passing to a subsequence, we can assume $(X_i,d_i,m_{G_i},p_i)$ converges in the pmGH-sense to $(X,d,m^*,p)$.\ By (2), $m^*=c\mathcal{H}^n$ for some $c>0$ and $(X,d,\mathcal{H}^n)$ is an RCD$(K,n)$ space.\ Also, $m$ and $\mathcal{H}^n$ are mutually absolutely continuous due to \cite{kell2017transport}.
		
		(4) By (1), $X$ is a compact connected nilpotent Lie group and hence, a torus.\ Since the metric is invariant and Riemannian, $X$ is a flat torus.
	\end{proof}
	
	\begin{remark}
		One may prove (4) without using the nilpotency obtained in (1).\ In fact, if $X$ is compact, then $G_i$ are finite groups and the orbits $G_i\cdot p_i$ are finite homogeneous metric spaces.\ Notice that $G_i\cdot p_i$ converges in the GH-sense to $X$.\ Then by \cite[Theorem 1.1]{gelander2013limits} and \cite[Theorem 2.2.4]{benjamini2017scaling}, $X$ is a torus with an invariant metric.\ Since $X$ is a Riemannian manifold, $X$ is a flat torus.
	\end{remark}
	
	Below, we will present a number of results that follow from Theorem \ref{thm-almost homogeneous rcd}.\ The following is readily obtained from Theorem \ref{thm-almost homogeneous rcd} and Theorem \ref{thm-zamora} (see also Remark \ref{rem-simply connected}).
	
	\begin{proposition}\label{prop-simply connected}
		Let $K\in\mathbb{R}$, $N\in(1,\infty)$ and $(X_i,d_i,m_i,p_i)$ be a sequence of almost homogeneous pointed RCD$(K,N)$ spaces converging to $(X,d,m,p)$ in the pmGH-sense.\ Assume that $\pi_1(X_i)$ are finite groups.\ Then $X$ is isometric to a simply connected nilpotent Lie group with a left invariant Riemannian metric.\ In particular, $X$ is diffeomorphic to $\mathbb{R}^n$, where $n=\dim(X)$.
	\end{proposition}
	
	If we additionally assume that $X_i$ are compact in Proposition \ref{prop-simply connected}, then the sequence $X_i$ must be collapsed.
	\begin{proposition}\label{prop-collapse}
		Let $(X_i,d_i,m_i,p_i)$ and $(X,d,m,p)$ be as in Proposition \ref{prop-simply connected}.\ Additionally, assume $X_i$ are compact and not single points.\ Then $\dim(X) \leq \liminf\limits_{i\to\infty} \dim_{ess}(X_i)-1 $.
	\end{proposition}
	
	\begin{proof}
		We can assume without loss of generality that $\dim_{ess}(X_i)=k \geq1$ for all $i$.\ Then by \cite[Theorem 1.5]{kitabeppu2019sufficient}, $\dim(X)\leq k$.\ Let $(X_i,d_i)$ be $(\eps_i,G_i)$-homogeneous with $\eps_i\to0$.\ Since $G_i$ are finite groups, we can apply \cite[Theorem A]{santos2020invariant} to obtain $G_i$-invariant measures $m_{G_i}$ such that $(X_i,d_i,m_{G_i})$ are also RCD$(K,N)$.\ Recall that the essential dimension is invariant under changes of measure (see \cite[Remark 2.12]{brena2023weakly} for instance).
		
		Let us now argue by contradiction.\ Suppose that $\dim(X)=k$.\ Then it follows from Theorem \ref{nss} that any sequence of small subgroups $W_i\leq\op{Iso}(X_i)$ will be trivial for large enough $i$.\ Then by Lemma \ref{lem:almost-homogeoeus-subgroups} and Lemma \ref{lem:free}, there exist subgroups $G_i^{\prime}\leq G_i$ such that $\diam(X_i/G_i^{\prime})\to0$ and $G_i^{\prime}$ acts freely on $X_i$ for large $i$.\ So $X_i\to X_i/G_i^{\prime}$ is a covering.\ Let $\tilde{X}_i$ be the universal cover of $X_i$ and $X_i/G_i^{\prime}$.\ Note that $\tilde{X}_i$ are compact.
		
		By \cite[Theorem 7.24]{galaz2018quotients}, $X_i/G_i^{\prime}$ endowed with the quotient metric and quotient measure is an RCD$(K,N)$ space.\ Then it follows from \cite[Theorem 2]{santos2023fundamental} that $\diam(\tilde{X}_i)\to0$.\ Hence $\diam(X_i)\to0$ and this leads to a contradiction.
	\end{proof}
	
	For an RCD$(K,N)$ space with $K>0$, the Bonnet-Myers theorem on RCD spaces \cite{sturm2006geometry2} will lead to a uniform diameter upper bound and the finiteness of the fundamental group.\ So the following corollary follows from Proposition \ref{prop-simply connected}.
	
	\begin{corollary}\label{cor-bonnet myers}
		Let $(X_i,d_i,m_i)$ be a sequence of almost homogeneous RCD$(K,N)$ spaces for some $K>0$ and $N\in(1,\infty)$.\ Then $\diam(X_i)\to0$.
	\end{corollary}
	\begin{remark}
		Corollary \ref{cor-bonnet myers} can also be obtained from (3) and (4) in Theorem \ref{thm-almost homogeneous rcd}.
	\end{remark}

	Recall that if an RCD$(0,N)$ space $(X,d,m)$ admits a discrete cocompact group $G\leq\op{Iso}(X)$, then $X$ splits as $\mathbb{R}^k\times Y$ where $Y$ is a compact RCD$(0,N-k)$ space (see Theorem \ref{gigli} and Remark \ref{rem-splitting}).\ We can derive the following corollary directly from Proposition \ref{prop-simply connected} and proof by contradiction.
	
	\begin{corollary}\label{cor-pan rong}
		Let $(X,d,m)$ be an RCD$(0,N)$ space for some $N\geq1$.\ Assume that $\pi_1(X)$ is finite and $G$ is a discrete subgroup of $\op{Iso}(X)$ with $\diam(X/G)<\infty$.\ Then $(X,d,m)$ is isomorphic to $(\mathbb{R}^k\times Y, d_{\mathbb{R}^k} \times d_Y, \mathcal{H}^k \otimes \nu)$, where $(Y,d_Y,\nu)$ is a compact RCD$(0,N-k)$ space with $\op{diam}(Y)\leq C(N)\cdot\op{diam}(X/G)$. 
	\end{corollary}
	
	\begin{proof}
		Without loss of generality, we can assume that $\diam(X/G)=1$.\ Let us argue by contradiction.\ Suppose that there is a sequence of $(1,G_i)$-homogeneous RCD$(0,N)$ spaces $(X_i,d_i,m_i)$ which are isomorphic to $\mathbb{R}^{k_i}\times Y_i$, where $Y_i$ are compact RCD$(0,N-k_i)$ spaces with finite fundamental groups and $\diam(Y_i)\to\infty$.\ We may also assume $k_i\equiv k$.
		
		Let $r_i=\diam(Y_i)$.\ Then $r_i^{-1}(\mathbb{R}^k\times Y_i)\to \R^k\times Y$ in the pGH-sense for some space $Y$ with $\diam(Y)=1$.\ Due to Proposition \ref{prop-simply connected}, this is a contradiction.
	\end{proof}
	
	The above corollary is also obtained by Pan-Rong in \cite{pan2018ricci}, where they considered Riemannian manifolds with $\Ric\geq0$.
	
	Theorem \ref{thm-almost nonnegative ricci-limit} is a direct consequence of Theorem \ref{thm-almost homogeneous rcd} and Theorem \ref{gigli}.

	\begin{proof}[Proof of Theorem \ref{thm-almost nonnegative ricci-limit}]
		It follows from Theorem \ref{thm-almost homogeneous rcd} and Theorem \ref{gigli} that $X$ is isometric to $\mathbb{R}^k\times N^{n-k}$ for some $0\leq k\leq n\leq N$, where $N^{n-k}$ is a nilpotent Lie group with a left invariant Riemannian metric which contains no lines.\ From the proof of Theorem \ref{thm-almost homogeneous rcd}, we know that $X$ is homogeneous which implies that $N^{n-k}$ is homogeneous.\ Then $N^{n-k}$ is compact and hence, is a torus.\ Since the metric is Riemannian and invariant, $N^{n-k}$ is a flat torus $\mathbb{T}^{n-k}$.
	\end{proof}
	
	\begin{remark}
		We note that Theorem \ref{thm-almost nonnegative ricci-limit} can also be derived from Theorem \ref{thm-almost homogeneous rcd} and \cite[Theorem 1.1]{huang2020non}.
	\end{remark}
	
	Theorem \ref{thm-almost nonnegative ricci-limit} and Proposition \ref{prop-simply connected} simply implies the following corollary.
	
	\begin{corollary}\label{cor-almost ricci positive}
		Let $(X_i,d_i,m_i,p_i)$ be a sequence of almost homogeneous RCD$(-\delta_i,N)$ spaces converging to $(X,d,m,p)$ in the pmGH-sense with $\delta_i\to0$.\ Suppose that $\pi_1(X_i)$ are finite groups.\ Then $X$ is isometric to $\mathbb{R}^{n}$ for some $n\leq N$.
	\end{corollary}
	
	
	
	The proof of Corollary \ref{cor-homeomorphic to tori} and Theorem \ref{thm-covering and almost homogeneous} is immediate.
	
	\begin{proof}[Proof of Corollary \ref{cor-homeomorphic to tori}]
		We argue by contradiction.\ Suppose that there is a sequence of almost homogeneous RCD$(K,N)$ spaces $(X_i,d_i,\mathcal{H}^N)$ with $\mathcal{H}^N(X_i)\geq v$ and $\diam(X_i)\leq D$ and all $X_i$ are not bi-H\"older homeomorphic to flat torus $\mathbb{T}^N$.\ Up to a subsequence, $X_i$ converges in the GH-sense to $Y$ and by Theorem \ref{thm-almost homogeneous rcd} (4) and Theorem \ref{thm:volume-continuity}, $Y$ is isometric to a flat torus $\mathbb{T}^N$.\ Then by Theorem \ref{thm:reifenberg-weak}, $X_i$ is bi-H\"older homeomorphic to $\mathbb{T}^N$ for sufficiently large $i$, which leads to a contradiction.
		
		Moreover, if $X$ is a Riemannian manifold, then we can show that $X$ is diffeomorphic to $\mathbb{T}^N$ by the same argument applying \cite[Theorem A.1.12]{cheeger1997structure} instead of Theorem \ref{thm:reifenberg-weak}.
	\end{proof}
	
	\begin{proof}[Proof of Theorem \ref{thm-covering and almost homogeneous}]
		Using Theorem \ref{thm-almost homogeneous rcd} (4), Theorem \ref{wang-fibration} and proof by contradiction, the proof is easily obtained.\ For the smooth case, substitute Theorem \ref{thm-local covering geometry} for  Theorem \ref{wang-fibration} and the conclusion follows from the same argument.
	\end{proof}


	\section{Topological rigidity of almost homogeneous non-collapsed RCD spaces}\label{sec-4}
	
	The main purpose of this section is to prove Theorem \ref{thm-rcd-maximal rank}.\ The proof will be divided into two parts (Theorem \ref{thm-rcd rigidity} and Theorem \ref{thm-rcd biHolder}), based on Zamora-Zhu's results in \cite{zamora2024topological} and Wang's arguments in \cite{wang2024non} respectively.\ In addition, we will simultaneously obtain a proof of Theorem \ref{thm-orbifold almost flat}.
	
	Let us first review the following theorem from \cite{zamora2024topological}.
	
	\begin{theorem}[\citeonline{zamora2024topological}]\label{thm-zamora zhu}
		For each $K\in\mathbb{R}$ and $N\geq1$, there is $\eps=\eps(K,N)$ such that if $(X,d,m)$ is an $(\eps,G)$-homogeneous RCD$(K,N)$ space, then $\op{rank}(G)\leq N$ and in the case of equality, $X$ is homeomorphic to $\mathbb{R}^N$.
	\end{theorem}
	The contractibility of $X$ is particular powerful when paired with the following theorem, which is an observation by Kapovitch in \cite{kapovitch2021mixed}.\ We refer to \cite{kapovitch2021mixed} (see also \cite{zamora2024topological}) for the proof.
	
	\begin{theorem}[\citeonline{kapovitch2021mixed}]\label{thm-Borel conj}
		Let $M$ be a closed aspherical topological manifold with $\pi_1(M)$ virtually nilpotent.\ Then $M$ is homeomorphic to an infranilmanifold.
	\end{theorem}
	
	By the definition of $\op{rank}(G)$ (see Definition \ref{def-rank}), the group $G$ in Theorem \ref{thm-zamora zhu} is virtually polycyclic.\ Indeed, due to Breuillard-Green-Tao's result \cite{breuillard2012structure}, $G$ is virtually nilpotent.
	
	 \begin{lemma}\label{lem-nilpotent group}
	     For each $K\in\mathbb{R}$ and $N\geq1$, there is $\eps=\eps(K,N)$ such that if $(X,d,m)$ is an $(\eps,G)$-homogeneous RCD$(K,N)$ space, then $G$ is finitely generated virtually nilpotent with $\op{rank}(G)\leq N$.
	 \end{lemma}
	 
	 \begin{proof}
	 	Fix $p\in X$.\ By \cite[Lemma 2.5]{zamora2024limits}, $G$ is generated by the set 
	 	\[S:=\left\lbrace g\in G :\ d(gp,p)\leq 3\cdot\diam(X/G)\right\rbrace. \]
	 	Then by \cite[Corollary 1.15]{breuillard2012structure}, there is small $\eps=\eps(K,N)$ such that $G$ is finitely generated virtually nilpotent.\ It follows from Theorem \ref{thm-zamora zhu} that $\op{rank}(G)\leq N$.
	 \end{proof}
	 
	 When $\rank(G)$ attains its maximum value $N$, $X$ is homeomorphic to $\mathbb{R}^N$ (Theorem \ref{thm-zamora zhu}) and in fact, the converse also holds.\ In addition, we can prove the first part of Theorem \ref{thm-rcd-maximal rank}, which demonstrate a set of conditions equivalent to maximal rank.\ Also, note that the first statement in Theorem \ref{thm-orbifold almost flat} is just a corollary.
	 
	 \begin{theorem}\label{thm-rcd rigidity}
	 	For each $K\in\mathbb{R}$ and $N\geq1$, there is $\epsilon=\epsilon(K,N), v=v(K,N)$ such that for any $(\eps,G)$-homogeneous RCD($K,N$) space $(X,d,m)$, the followings are equivalent:
	 	\begin{enumerate}
	 		\item $X$ is homeomorphic to $\mathbb{R}^N$;
	 		\item $X$ is a contractible topological $N$-manifold without boundary;
	 		\item $\op{rank}(G)$ is equal to $N$;
	 		\item $X$ is simply connected and $\mathcal{H}^N(B_1(x))\geq v$ for some $x\in X$;
	 		\item $\pi_1(X)$ is finite and $\mathcal{H}^N(B_1(x))\geq v$ for some $x\in X$.
	 	\end{enumerate}
	 \end{theorem}
	 
	 \begin{proof}
	 	We will prove that (3)$\Rightarrow$(4)$\Rightarrow$(5)$\Rightarrow$(3) and (1)$\Rightarrow$(2)$\Rightarrow$(3)$\Rightarrow$(1).
	 	
	 	(3)$\Rightarrow$(4): By Theorem \ref{thm-zamora zhu} and Theorem \ref{thm:dim}, $X$ is simply connected and $m=c\mathcal{H}^N$ for some $c>0$.\ We only need to show $\mathcal{H}^N(B_1(x))\geq v$ for some $x\in X$.
	 	
	 	By contradiction, we assume that there is a sequence of $(\eps_i,G_i)$-homogeneous RCD$(K,N)$ spaces $(X_i,d_i,\mathcal{H}^N)$ with $\eps_i\to0$, such that $\rank(G_i)=N$ and $\mathcal{H}^N(B_1(x_i))\to0$ for some sequence $x_i\in X_i$.\ By compactness, Theorem \ref{thm-almost homogeneous rcd} and Theorem \ref{thm:volume-continuity}, we can assume $(X_i,x_i)$ converges in the pGH-sense to $(X,x)$, where $X$ is a nilpotent Lie group of dimension $n\leq N-1$.\ Let $G_i^{\prime}$ be as in Lemma \ref{lem:almost-homogeoeus-subgroups}.\ Then by Remark \ref{rem:rank-inequality}, $\rank(G_i^{\prime})\leq N-1$ which implies $\rank(G_i)\leq N-1$.\ This leads to a contradiction.
	 	
	 	(4)$\Rightarrow$(5): This is trivial.
	 	
	 	(5)$\Rightarrow$(3): Notice that $\dim_{\mathcal{H}}(X)=N$ and by Theorem \ref{thm:dim}, $m=c\mathcal{H}^N$ for some $c>0$.\ We then argue by contradiction.\ Assume that there is a sequence of $(\eps_i,G_i)$-homogeneous RCD$(K,N)$ spaces $(X_i,d_i,\mathcal{H}^N)$ with $\eps_i\to0$, such that $\pi_1(X_i)$ are finite, $\rank(G_i)<N$ and $\mathcal{H}^N(B_1(x_i))\geq v>0$ for some sequence $x_i\in X_i$.\ Due to compactness, Theorem \ref{thm-almost homogeneous rcd} and Theorem \ref{thm:volume-continuity}, we can assume $(X_i,x_i)$ converges in the pGH-sense to $(X,x)$, where $X$ is a nilpotent Lie group of dimension $N$.\ Then by Lemma \ref{lem:malcev-construction} and Remark \ref{rem-simply connected}, $\rank(G_i)=N$ which leads to a contradiction.
	 	
	 	(1)$\Rightarrow$(2): This is trivial.
	 	
	 	(2)$\Rightarrow$(3): By Lemma \ref{lem-nilpotent group}, $G$ is a finitely generated virtually nilpotent group.\ Then $G$ contains a torsion free nilpotent subgroup of finite index (see \cite{kargapolov1979fundamentals}), denoted by $\Gamma$.\ Notice that $\Gamma\leq\op{Iso}(X)$ is a discrete group acting freely on $X$.\ Since $X$ is a contractible topological $N$-manifold, $X/\Gamma$ is a closed aspherical topological manifold.\ Then by Theorem \ref{thm-Borel conj}, $X/\Gamma$ is homeomorphic to an $N$-dimensional nilmanifold.\ So $\Gamma=\pi_1(X/\Gamma)$ has rank $N$ and thus, $\rank(G)=N$.
	 	
	 	(3)$\Rightarrow$(1): This follows from Theorem \ref{thm-zamora zhu}.
	 \end{proof}
	
	
	We now proceed to prove the last statement in Theorem \ref{thm-rcd-maximal rank} and Theorem \ref{thm-orbifold almost flat}.\ Notice that we only need to work on non-collapsed RCD$(K,N)$ spaces $(X,d,\mathcal{H}^N)$.\ Also, if $G$ is isomorphic to an almost-crystallographic group of dimension $N$, then $\op{rank}(G)=N$.\ Therefore, we only need to show the following theorem.
	
	\begin{theorem}\label{thm-rcd biHolder}
		For each $K\in\mathbb{R}$ and $N\geq1$, there is $\epsilon=\epsilon(K,N)$ such that for any $(\eps,G)$-homogeneous RCD($K,N$) space $(X,d,\mathcal{H}^N)$, if $G$ does not contain a non-trivial finite normal subgroup and $\op{rank}(G)=N$, then $X/G$ is bi-H\"older homeomorphic to an $N$-dimensional infranil orbifold $\mathcal{N}/\Gamma$, where $\mathcal{N}$ is a simply connected nilpotent Lie group endowed with a left invariant metric and $G$ is isomorphic to $\Gamma$.
		
		Furthermore, if $X$ is a smooth Riemannian manifold, then the Riemannian orbifold $X/G$ is diffeomorphic to an $N$-dimensional infranil orbifold.
	\end{theorem}
	
	First note that Theorem \ref{thm-almost crystallographic} and Lemma \ref{lem-nilpotent group} will imply that for a small $\eps$, the group $G$ in Theorem \ref{thm-rcd biHolder} is isomorphic to an almost-crystallographic group of dimension $N$.
	
	Assume that Theorem \ref{thm-rcd biHolder} does not hold.\ Then after rescaling on metrics, there is a sequence of $(\eps_i,G_i)$-homogeneous RCD$(-\eps_i,N)$ spaces $(X_i,d_i,\mathcal{H}^N)$ with $\eps_i\to0$, such that any $G_i$ is isomorphic to an almost-crystallographic group of dimension $N$ and $X_i/G_i$ is not bi-H\"older homeomorphic to any infranil orbifold of the form $\nn_i/G_i$.
	
	Let $G_i^{\prime}$ be the bounded index normal subgroups of $G_i$ in Lemma \ref{lem:almost-homogeoeus-subgroups}.\ Then $\op{rank}(G_i^{\prime})=N$ and it follows from \cite[Lemma 2.6]{zamora2024limits} that $\diam(X_i/G_i^{\prime})\to0$.\ Due to Remark \ref{rem:rank-inequality} and Corollary \ref{cor-almost ricci positive}, we have the following diagram:
	    \[\xymatrix{
			(X_{i},p_{i},G_{i}^{\prime}) \ar[r]^{eqGH} \ar[d] & (\R^N,0,G) \ar[d] \\
			X_{i}/G_{i}^{\prime} \ar[r]^{GH}
			&  \text{ pt }.} \]
    By Lemma \ref{lem:almost-translations}, $G$ acts freely and transitively on $\R^N$ and hence, $G=\R^N$.
    
    By Theorem \ref{nss}, the groups $G_i$ admit on non-trivial small subgroups.\ Then by Lemma \ref{lem:malcev-construction} and Remark \ref{rem:rank-inequality}, $G_i^{\prime}$ are torsion free nilpotent groups.\ Let $\nn_i$ be the Mal'cev completion of $G_i^{\prime}$.\ It follows from Proposition \ref{prop-Malcev} that $G_i$ is an almost-crystallographic group modeled on $\nn_i$ for each $i$.


	Our goal is to find a left invariant metric on $\mathcal{N}_i$ so that $X_i/G_i$ is bi-H\"older homeomorphic to $\nn_i/G_i$.\ The proof is essentially the same as in \cite{wang2024non}.\ For the convenience of readers, we give the construction of the left invariant metric on $\mathcal{N}_i$ in \cite[Lemma 4.5]{wang2024non}.
	
	\begin{lemma}[\citeonline{wang2024non}]\label{metric}
		Let $(X_i,d_i,p_i,G_i^{\prime})$ and $\nn_i$ be as above.\ For any $\epsilon\in (0,1)$ and large $i$, $\mathcal{N}_i$ admits a left invariant metric $g_{\mathcal{N}_i}$ with  $\mathrm{inj}_{\mathcal{N}_i} \ge \frac{1}{\epsilon}$.\ Moreover, there is $\epsilon_i \to 0$ so that $\forall g \in \mathcal{N}_i$, $B_{\frac{1}{\epsilon}}(g) \subset \mathcal{N}_i$ is $\epsilon_i$-$C^4$-close, by $\mathrm{exp}_g^{-1}$, to the $\frac{1}{\epsilon}$-ball in $T_{g}\mathcal{N}_i$ with the flat metric.
	\end{lemma}
	
	\begin{proof}
	Note that the groups $G_i^{\prime}$ admit no non-trivial small subgroups and $\op{rank}(G_i^{\prime})=N$.\ Then by Lemma \ref{lem:malcev-construction} and Remark \ref{rem:rank-inequality}, for $i$ large enough there are generators $u_{1,i}, \ldots , u_{N,i} \in G_ i ^{\prime }$, and $C_{1,i}, \ldots , C_{N,i} \in \mathbb{R}^+$ with the following properties:
	\begin{enumerate}
		\item There are polynomials $Q_i : \mathbb {R}^N \times \mathbb{R}^N\to \mathbb{R}^N$ of degree $\leq d(N)$ giving the group structures on $\R^N$ by $x_1 \cdot x_2 = Q_i (x_1, x_2)$ such that for each $i$, $G_ i ^{\prime}$ is isomorphic to the group $(\Z^N,Q_i|_{\Z^N\times\Z^N})$ and the group $(\R^N,Q_i)$ is isomorphic to $\nn_i$.
		\item There is $C > 0 $ such that the set
		\[  
		P_i : = P (  u_{1,i}, \ldots , u_{N,i} ; C_{1,i}, \ldots , C_{N,i} )  \subset  G_ i ^{\prime} 
		\]  
		is a nilprogression in $C$-normal form with $\thi (P_i ) \to \infty $.  \label{property:zamora-3}
		\item  \label{property:zamora-4} For each $\varepsilon > 0 $ there is $\delta > 0 $ such that 
		\begin{align}
			G( \delta  P_ i  )    \subset  \{ & g \in G_i ^{\prime} \, \vert \, d_i(g p_i, p_i ) \leq \varepsilon \} ,  \label{eq:grid-small} 
		\end{align}
		for $i$ large enough.    
	\end{enumerate}
	
	
	By \eqref{eq:grid-small}, there is $\delta_1 > 0 $ with $G(\delta_1 P_i) \subset \{ g \in G_ i ^{\prime} \vert d_i(gp_i,p_i) \leq 1 \} $.\ Hence there is an integer $D \in \mathbb{N} $ so that 
	\begin{equation} \label{eq:grid-small-2}
		G(P_i) \subset  G(\delta_1  P_i) ^{D} \subset   \{ g \in G_ i ^{\prime} \vert d_i(gp_i, p_i) \leq D \} .   
	\end{equation}
	Let $g_{j,i}:=u_{j,i}^{\lfloor \frac{C_{j,i}}{C}\rfloor }$ and 
	$v_{j,i} : = \log (g_{j,i})\in T_e \nn_i $.\ Then $\{ v_{1,i}, \ldots , v_{N,i}\} $ is a strong Mal'cev basis of the Lie algebra $T_{e}\nn_i$ (see \cite{zamora2024limits}).
	
	Notice that the groups $G_i^{\prime}$ converge equivariantly to the group of translations in $\mathbb{R}^N$.\ By \eqref{eq:grid-small-2}, after passing to a subsequence, for each $j \in \{ 1, \ldots , N \}$ we can assume $g_{j,i}$ converges equivariantly to some $v_j \in \mathbb{R}^N$.\ We may identify $\R^N$ with its Lie algebra and $\{ v_1, \ldots , v_N \}$ is a basis of $\mathbb{R}^N$.
	
	Define the left invariant metric $g_{\mathcal{N}_i}$ by the inner product on $T_e\nn_i$ as following:
	$$g_{\mathcal{N}_i}(v_{j_1,i},v_{j_2,i}) = \langle v_{j_1}, v_{j_2} \rangle,$$ 
	where $1 \le j_1,j_2 \le N$ and the right-hand side is the inner product in $\mathbb{R}^N$.

	Since $\{ v_{j,i} , 1 \le j \leq N\}$ is a strong Malcev basis of $T_e\nn_i$, for any $1 \le j_1 < j_2 \le N$, 
	\begin{align}
		[v_{j_1,i},v_{j_2,i}]= \sum_{j=j_2+1}^n a_{j_1j_2,i}^j\ v_{j,i}.
	\end{align}
	It is proven in \cite[Lemma 2.64 and Proposition 8.2]{zamora2024limits} that the structure coefficients of $T_e\nn_i$ with respect to the basis $\{ v_{1,i}, \ldots , v_{N,i}\} $ converge to the structure coefficients in $\mathbb{R}^N$ with respect to $\{ v_1, \ldots , v_N \}$ as $i \to \infty$.\ Thus $a_{j_1j_2,i}^j \to 0$ as $i \to \infty$, since the limit group is abelian.\ Define $a_{j_1j_2,i}^j=0$ if $j \le j_1$ or $j \le j_2$.\ Then for any $ 1 \le j_1,j_2,j_3 \le N$,
	$$g_{\mathcal{N}_i}(\nabla_{v_{j_1,i}} v_{j_2,i} , v_{j_3,i}) = \frac{1}{2} (a_{j_1j_2,i}^{j_3} - a_{j_2j_3,i}^{j_1} + a_{j_3j_1,i}^{j_2}).$$
		
	Note that all terms on the right-hand side are constant (depending on $i$) and converge to $0$ as $i \to \infty$.\ In particular, the covariant derivatives of the Riemannian curvature tensor $g_{\mathcal{N}_i}$ satisfy 
	$$|(\nabla^{g_{\mathcal{N}_i}})^k Rm_{g_{\mathcal{N}_i}}| \le \epsilon_i, \ 0 \le k \le 3,$$
	where $\epsilon_i \to 0$.\ Then one can easily verify that this metric fulfills the conditions.
	\end{proof}
	
	\begin{remark}
		In \cite{wang2024non}, Wang used the results in \cite{breuillard2012structure,zamora2024limits} to embed $G_i'$ as a lattice in a simply connected nilpotent Lie group $\nn_i$.\ In fact, since $G_i'$ is torsion free nilpotent, one can directly use its Mal'cev completion.\ One can see Remark \ref{cor-rcd} for the reason why the groups $G_i$ and $G_i'$ are torsion free in Wang's theorem.
	\end{remark}
	
	From now on, $\mathcal{N}_i$ is always endowed with the metric $g_{\mathcal{N}_i}$ constructed in Lemma \ref{metric}.
	
	Define $G_i^{\prime}(p_i,D):=\{g \in G_i ^{\prime}\vert d_i(gp_i, p_i)\leq D\}$.\ Assume that a pseudo-group $G$ acts on two metric spaces $X_1,X_2$ separately by isometries.\ Following \cite{wang2024non}, we say a map $h:X_1 \to X_2$ is $\epsilon$-almost $G$-equivariant if $d(h(gx),gh(x)) < \epsilon$ for any $ x \in X_1, g \in G$.
	
	The following two lemmas come from \cite{wang2024non,wang2023limit}.
	
	\begin{lemma}[\citeonline{wang2024non}]\label{localeg}
		For any $ \epsilon > 0$, let $B_{\frac{1}{\epsilon}}(p_i) \subset X_i$ and $B_{\frac{1}{\epsilon}}(e) \subset \mathcal{N}_i$.\ Then there exists an $\epsilon_i$-GHA $h_i': B_{\frac{1}{\epsilon}}(p_i) \to B_{\frac{1}{\epsilon}}(e)$ which is $\epsilon_i$-almost $G_i^{\prime}(p_i,{\frac{1}{\epsilon}})$-equivariant if it is well-defined, where $\epsilon_i \to 0$ as $i \to \infty$.
	\end{lemma}
	
	
	\begin{lemma}[\citeonline{wang2023limit,wang2024non}]\label{extension}
		The map $h_i^{\prime}$ in Lemma \ref{localeg} can be extended to a global map $h_i:X_i \to \mathcal{N}_i$, which is an $\epsilon_i$-GHA on any $\frac{1}{\epsilon}$-ball and $\epsilon_i$-almost $G_i'$-equivariant with $\eps_i\to0$.
	\end{lemma}

	Now we can follow the arguments in \cite{wang2024non} to complete the proof of Theorem \ref{thm-rcd biHolder}.
	\begin{proof}[Proof of Theorem \ref{thm-rcd biHolder}] 
		Let us argue by contradiction.\ There is a sequence of $(\eps_i,G_i)$-homogeneous RCD$(-\eps_i,N)$ spaces $(X_i,d_i,\mathcal{H}^N)$ with $\eps_i\to0$ and $\op{rank}(G_i)=N$.\ Also, we have already established the following diagram:
		\[\xymatrix{
			(X_{i},p_{i},G_{i}^{\prime}) \ar[r]^{eqGH} \ar[d] & (\R^N,0,\R^N) \ar[d] \\
			X_{i}/G_{i}^{\prime} \ar[r]^{GH}
			&  \text{ pt },} \]
		where for each $i$, $G_i'$ is a normal subgroup of bounded index in $G_i$, embedded as a lattice in an $N$-dimensional simply connected nilpotent Lie group $\mathcal{N}_i$ and $G_i$ is an almost-crystallographic group modeled on $\nn_i$.\ We also assumed that none of $X_i/G_i$ is bi-H\"{o}lder homeomorphic to the infranil orbifold $\nn_i/G_i$. 
		
		By the construction of the metric $g_{\nn_i}$ in Lemma \ref{metric}, the lattice $G_i'$ is $\epsilon_i$-dense in $\mathcal{N}_i$.\ Since $\diam(X_i/G_i')\to0$, the map $h_i$ in Lemma \ref{extension} is also $\epsilon_i$-almost $G_i$-equivariant.
		
		By the same arguments in the proof of \cite[Theorem A]{wang2024non}, we can assume that $G_i$ acts on $\mathcal{N}_i$ by isometries and for any small $\eps>0$, there is a normal subgroup $G_i''$ in $G_i'$ of finite index, which is also normal in $G_i$, so that $G_i'' \cap B_{\frac{1}{\epsilon}}(e) = \{e\}$.
		
		Since $G_i'' \cap B_{\frac{1}{\epsilon}}(e) = \{e\}$, we can apply Lemma \ref{metric} to conclude that the injective radius of  $\mathcal{N}_i/G_i''$ is at least ${\frac{1}{\epsilon}}$.\ For any $y \in \mathcal{N}_i/G_i''$, $B_{\frac{1}{\epsilon}}(y) \subset \mathcal{N}_i/G_i''$ is $\epsilon_i$-$C^4$-close to the $\frac{1}{\epsilon}$-ball in the tangent space $T_{y}(\mathcal{N}_i/G_i'')$ with the flat metric. 
		
		Since $h_i$ is $\epsilon_i$-almost $G_i$-equivariant, we can reduce $h_i$ to a map 
		$$\bar{h}_i: X_i/G_i''\to \mathcal{N}_i/G_i'',$$
		which is an $\epsilon_i$-GHA on any ${\frac{1}{\epsilon}}$-ball and $\epsilon_i$-almost $G_i/G_i''$-equivariant.
		
		Since $G_i/G_i''$ is finite, we can apply \cite[Theorem 3.5]{wang2024non} to 
		$$\bar{h}_i: (X_i/G_i'',G_i/G_i'')\longrightarrow (\mathcal{N}_i/G_i'',G_i/G_i'').$$ 
		Thus there is a $(G_i/G_i'')$-equivariant map $f_{G_i/G_i''}: X_i/G_i'' \to \mathcal{N}_i/G_i''$, which is harmonic $(N,\Phi(\epsilon|N))$-splitting on any $\frac{1}{5\epsilon}$-ball.\ Then by Theorem \ref{Rei},
		$$(1-\Phi(\epsilon|N))d_i(x,y)^{1+\Phi(\epsilon|N)} \le d(f_{G_i/G_i''}(x),f_{G_i/G_i''}(y)) \le (1+\Phi(\epsilon|N))d_i(x,y),$$
		for any $x,y \in X_i/G_i''$ with $d_i(x,y) \le \frac{1}{10\epsilon}$. 
		
		Since $f_{G_i/G_i''}$ is $(G_i/G_i'')$-equivariant, it can be reduced to a bi-H\"{o}lder map on the quotient space $f: X_i/G_i \to \mathcal{N}_i/G_i$.\ This leads to a contradiction to the assumption.
		
		Furthermore, if $X_i$ is a Riemannian manifold, then $X_i/G_i''$ is also a Riemannian manifold, since the group $G_i''$ acts freely on $X_i$.\ So by Theorem \ref{Rei}, the $(G_i/G_i'')$-equivariant map $f_{G_i/G_i''}: X_i/G_i'' \to \mathcal{N}_i/G_i''$ restricted on any $\frac{1}{10\epsilon}$-ball, is a diffeomorphism onto its image.\ Notice that $X_i/G_i''\to X_i/G_i$ and $\nn_i/G_i''\to\nn_i/G_i$ are orbifold coverings.\ Hence, the reduced map $f: X_i/G_i \to \mathcal{N}_i/G_i$ is a diffeomorphism between orbifolds.\ This completes the proof.
	\end{proof}
	
	Combining Theorem \ref{thm-rcd rigidity} and Theorem \ref{thm-rcd biHolder}, both Theorem \ref{thm-rcd-maximal rank} and Theorem \ref{thm-orbifold almost flat} are readily obtained.
	
	\begin{remark}
		In the last statement of Theorem \ref{thm-orbifold almost flat}, it is expected that the assumptions requiring a good orbifold and an orbifold fundamental group without non-trivial finite normal subgroups can be eliminated.\ Due to \cite[Proposition 1.4]{ding2011restriction}, any almost flat orbifold is an infranil orbifold.\ So it might be more natural to seek a nearby almost flat metric under the conditions in Theorem \ref{thm-orbifold almost flat}.\ This is achieved in the manifold case via Ricci flow smoothing techniques (see \cite{huang2020collapsed}).
	\end{remark}
	

	\begin{remark}\label{cor-rcd}
		Although in Theorem \ref{thm-rcd-maximal rank}, we assume that the group $G$ does not contain a non-trivial finite normal subgroup, Theorem \ref{thm-rcd-almostflat} is still a corollary of Theorem \ref{thm-rcd-maximal rank}.\ This is due to the fact that if $Y$ is a closed aspherical topological manifold, then $\pi_1(Y)$ is torsion free.\ Notice that any closed topological manifold is homotopy equivalent to a CW complex \cite{kirby1969triangulation} and the fundamental group of an aspherical finite-dimensional CW complex is torsion free \cite{luck2012aspherical}.
	\end{remark}


	\section{Rigidity and regularity of almost homogeneous Einstein metrics}\label{sec-5}
	
	In this section, we mainly focus on the rigidity and $\eps$-regularity for almost homogeneous Riemannian orbifolds and manifolds with bounded Ricci curvature.\ We first give the proof of Theorem \ref{thm-einstein orbifold}, which is an orbifold verion of \cite[Theorem 0.2]{si2024rigidity}.
	
	\begin{proof}[Proof of Theorem \ref{thm-einstein orbifold}]
		
		Since (1)$\Rightarrow$(2)$\Rightarrow$(3) is trivial and (3) is equivalent to (4) by Theorem \ref{thm-rcd-maximal rank}, it suffices to show that (3) and (4) together imply (1).
		
		Let us argue by contradiction.\ Suppose that there is a sequence of non-flat Einstein $n$-orbifolds $(\mathcal{O}_i,g_i)$ such that $\op{Ric}_{g_i}\equiv\lambda_i$ with $\lambda_i\ge -(n-1)$, $\op{diam}(\mathcal{O}_i,g_i)\to 0$, and satisfying (3) and (4) in Theorem \ref{thm-einstein orbifold}.
		
		Consider the universal orbifold covers $(\tilde{\Or}_i,\tilde{g}_i)$.\ Due to Corollary \ref{cor-bonnet myers} and Corollary \ref{cor-pan rong}, we can assume that $\lambda_i< 0$.\ Then up to a rescaling on metrics, we can further assume that $\lambda_i=-(n-1)$.\ Note that $\orb_i$ still converges to a point and $\op{rank}(\pi_1^{orb}(\orb_i))=n$.\ By Theorem \ref{thm-rcd-maximal rank}, $\op{vol}_{\tilde{g}_i}(B_1(\tilde{x}_i))\geq v^{\prime}(n)>0$ for some $\tilde{x}_i\in\tilde{\Or}_i$.\ Up to a subsequence, we have the following pmGH-convergence by Theorem \ref{thm:volume-continuity},
		\[(\tilde{\Or}_i,\tilde{g}_i,\op{vol}_{\tilde{g}_i},\tilde{x}_i) \xrightarrow{pmGH} (\tilde{X},\tilde{g},\mathcal{H}^n,\tilde{x}),\]
		where by Proposition \ref{prop-simply connected}, $\tilde X $ is isometric to a simply connected nilpotent Lie group with a left invariant Riemannian metric, denoted by $\tilde{g}$.
		
		On the other hand, for any $i$, $(|\tilde \orb_i|_{reg},g_{i, reg})$ is a smooth open Riemannian manifold with $\Ric\equiv-(n-1)$ and $\op{vol}_{\tilde{g}_i}(B_r(\tilde y_i)\cap|\tilde \orb_i|_{reg}) = \op{vol}_{\tilde{g}_i}(B_r(\tilde y_i))$ for any $\tilde{y}_i\in \tilde{\Or}_i$ and $r>0$.\ By the standard Schauder estimate, $(|\tilde \orb_i|_{reg},g_{i, reg})$ converges in the $C^{\infty}_{loc}$-norm to a full measure subset of $(\tilde{X},\tilde{g})$ (see \cite{cheeger1997structure}).\ Since $\tilde{X}$ is a Riemannian manifold, $\Ric_{\tilde{X}}\equiv-(n-1)$.
		
		By \cite[Theorem 2.4]{milnor1976curvatures}, any left invariant Riemannian metric of a nilpotent but not abelian Lie group has both directions of strictly negative and positve Ricci curvature.\ Thus, $(\tilde X, \tilde g)$ must be isometric to $\mathbb{R}^n$, which leads to a contradiction.
	\end{proof}
	
	Recall that if $(M,g)$ is a Riemannian manifold and $G$ is a discrete subgroup of $\op{Iso}(M)$, then $M/G$ admits a natural orbifold structure.\ Let $\tilde{M}$ be the universal cover of $M$.\ Then $\tilde{M}$ is also the universal orbifold cover of the good orbifold $M/G$.\ Moreover, if $M$ is simply connected, then $\pi_1^{orb}(M/G)=G$.\ Therefore, the following corollary is readily derived from Theorem \ref{thm-rcd-maximal rank} and Theorem \ref{thm-einstein orbifold}.
	
	\begin{corollary}\label{cor-rigidity}
		There is $\epsilon=\epsilon(n)>0,v=v(n)>0$ such that if an $(\epsilon,G)$-homogeneous Einstein $n$-manifold $(M,g)$ satisfies $\op{Ric}_g=\lambda g$ with $\lambda\geq -(n-1)$, then the followings are equivalent:
			\begin{enumerate}
						\item $\op{vol}(B_1(x))\geq v$ for some $x\in M$, and $M$ is simply connected;
						\item $\op{rank}(G)=n$;
						\item $M$ is diffeomorphic to $\mathbb{R}^n$;
						\item $M$ is isometric to $\mathbb{R}^n$.
			\end{enumerate}
		In particular, if we only assume $\op{vol}(B_1(x))\geq v$ for some $x\in M$, then $M$ is flat.
	\end{corollary}
	The following proposition is a quantitative rigidity version of Corollary \ref{cor-rigidity}.
	
	\begin{proposition}
		Given $v>0$ and $p\in(1,\infty)$, for any $\delta>0$, there is $\epsilon=\epsilon(n,v,p,\delta)$ such that if an $(\epsilon,G)$-homogeneous $n$-manifold $(M,g)$ satisfies $|\op{Ric}-\lambda g|\leq\eps$ with $\lambda\geq -(n-1)$ and $\op{vol}(B_1(x))\geq v$, then $\int_{B_1(x)}|Rm|^p\leq \delta$.
	\end{proposition}
	
	\begin{proof}
		Argue by contradiction.\ Suppose that there exists $\delta_0>0$ such that for any $\eps_i\to0$, there is a sequence of $(\eps_i,G_i)$-homogeneous pointed $n$-manifolds $(M_i,g_i,x_i)$ satisfying $|\op{Ric}-\lambda_i g_i|\leq\eps_i$ with $\lambda_i \geq -(n-1)$, $\op{vol_{g_i}}(B_1(x_i))\geq v$ and $\int_{B_1(x_i)}|Rm|^p > \delta_0$.
		
		By Corollary \ref{cor-bonnet myers}, we can assume that $\lambda_i$ converges to some $\lambda_\infty \leq0$.\ Up to a subsequence, we have the following pGH-convergence
		\[(M_i,g_i,x_i)\xrightarrow{pGH}(X,d,x),\]
		where $(X,d)$ is isometric to a Riemannian manifold with $\Ric\geq\lambda_\infty$ by Theorem \ref{thm-almost homogeneous rcd}.\ Note that $(M_i,g_i,x_i)$ converges in the pointed $C^{1,\alpha}\cap W^{2,q}$-topology to $(X,g_X,x)$, where the metric $g_X$ is a weak solution of the Einstein equation 
		$$\Delta g + Q(g,\partial g)=\lambda_\infty g,$$ 
		under harmonic coordinate charts (see \cite{anderson1990convergence,cheeger1997structure}).\ Hence, $g_X$ is a smooth metric and $X$ is an Einstein manifold with $\Ric_{g_X}=\lambda_\infty$.\ By the same arguments in the proof of Theorem \ref{thm-einstein orbifold}, $X$ is a flat manifold.\ Since $(M_i,g_i,x_i)$ converges in the pointed $C^{1,\alpha}\cap W^{2,q}$-topology to $(X,g_X,x)$ for any $0<\alpha<1$ and $1<q<\infty$, we have $\int_{B_1(x_i)}|Rm|^p \to 0$.\ This leads to a contradiction.
	\end{proof}
	
	Let $(M,g)$ be a Riemannian manifold.\ Recall that the $C^k$-harmonic radius at $x\in M$ is defined to be the largest $r>0$ such that there exists a harmonic coordinate system on $B_r(x)$ with $C^k$-control on the metric tensor.\ Harmonic coordinates have an abundancce of good properties when it comes to regularity issues.\ We refer to \cite{petersen1997convergence} for a nice introduction.\ In particular, if the Ricci curvature is uniformly bounded, then in harmonic coordinates, the metric $g_{ij}$ has \emph{a priori} $C^{1,\alpha}\cap W^{2,q}$-bounds for any $\alpha\in(0,1)$ and $q\in(1,\infty)$.
	
	Let us now proceed to prove the $\eps$-regularity theorem (Theorem \ref{thm-regularity}).
	
	\begin{proof}[Proof of Theorem \ref{thm-regularity}]
		Let us argue by contradiction.\ Suppose that there exists a sequence of almost homogeneous pointed $n$-orbifolds $(\orb_i,g_i,x_i)$ satisfying $|\Ric_{g_i}|\leq n-1$, $\op{vol}(B_1(x_i))\geq v$ and $\int_{B_1(x_i)}|Rm|^p\to\infty$.\ Then up to a subsequence, $(\orb_i,g_i,x_i)$ converges in the pGH-sense to $(X,d,p)$.\ By Theorem \ref{thm-almost homogeneous rcd}, $X$ is a Riemannian manifold.\ Hence, $B_1(x_i)\cap|\orb_i|_{reg}$ converges to an open subset of $(X,d,p)$ in the $C^{1,\alpha}\cap W^{2,q}$-norm for any $\alpha\in(0,1)$ and $q\in(1,\infty)$, which implies that $\int_{B_1(x_i)}|Rm|^p=\int_{B_1(x_i)\cap|\orb_i|_{reg}}|Rm|^p$ is bounded.\ This leads to a contradiction.
		
		For the last statement, notice that the $C^1$-harmonic radius $r_h$ is continuous under $C^{1,\alpha}$-topology and hence, a similar proof applies.
	\end{proof}
	
	The proof of Theorem \ref{thm-bounded Ricci almostflat} is readily obtained.
	
	\begin{proof}[Proof of Theorem \ref{thm-bounded Ricci almostflat}]
		By \cite[Corollary 1.2]{dai2000integral}, (1) implies (2).\ It follows from Theorem \ref{thm-almostflat} that (2), (3) and (4) are equivalent.\ Then Theorem \ref{thm-regularity} shows that (4) implies (1).\ One can adjust the constants to make (1), (2), (3) and (4) equivalent.
	\end{proof}

	\bibliographystyle{plain}
	\bibliography{almosthomogeneous}
\end{document}